\documentclass[11pt,draft]{amsart}
\usepackage{amssymb}
\usepackage[matrix,arrow]{xy}

\theoremstyle{plain}
\newtheorem*{assumption}{Assumption}
\numberwithin{equation}{section}

\theoremstyle{plain}
\newtheorem{bigtheorem}{Theorem}
\newtheorem{bigcorollary}[bigtheorem]{Corollary}
\newtheorem{bigproposition}[bigtheorem]{Proposition}
\newtheorem{thm}{Theorem}[section]
\newtheorem{theorem}[thm]{Theorem}
\newtheorem{cor}[thm]{Corollary}
\newtheorem{corollary}[thm]{Corollary}
\newtheorem{lemma}[thm]{Lemma}
\newtheorem{prop}[thm]{Proposition}
\newtheorem{proposition}[thm]{Proposition}

\theoremstyle{definition}
\newtheorem{defn}[thm]{Definition}
\newtheorem{definition}[thm]{Definition}
\newtheorem{rem}[thm]{Remark}
\newtheorem{remark}[thm]{Remark}

\newtheorem{convention}[thm]{Convention}

\author{Paul Seidel}

\newcommand{\standardspacing}{\parindent0cm \parskip1ex}

\newcommand{\R}{\mathbb{R}}
\newcommand{\Z}{\mathbb{Z}}

\newcommand{\C}{\mathbb{C}}
\newcommand{\N}{\mathbb{N}} 

\newcommand{\ind}{\mathrm{ind }}

\newcommand{\iso}{\cong}           

\newcommand{\smooth}{C^\infty}

\renewcommand{\o}{\omega}
\renewcommand{\O}{\Omega}
\newcommand{\mo}{(M,\o)}

\newcommand{\Aut}{\mathrm{Aut}}
\newcommand{\Ham}{\mathrm{Ham}}

\newcommand{\wplus}{$(W^+)$}
\newcommand{\loops}{\mathcal{L}M}
\newcommand{\tloops}{\widetilde{\loops}}
\newcommand{\quantum}{QH_*\mo}
\newcommand{\inv}{QH_*\mo^{\times}}
\renewcommand{\sec}{\mathcal{S}}
\newcommand{\ev}{\mathrm{ev}}
\newcommand{\novi}{\Lambda}
\newcommand{\hgroup}{\pi_2(M)}
\newcommand{\qp}{\ast}
\newcommand{\floer}{HF_*\mo}
\newcommand{\id}{\mathit{Id}}
\newcommand{\group}{\widetilde{G}}
\newcommand{\pair}{(g,\tilde{g})}
\newcommand{\pairofpants}{\qp_{\mathrm{PP}}}
\renewcommand{\H}{\mathcal H}
\newcommand{\glue}{\#}

\newcommand{\J}{\mathbf{J}}
\newcommand{\modJ}[1]{\mathcal{J}(#1;M,\o)}
\newcommand{\modJreg}[1]{\mathcal{J}^{\mathrm{reg}}(#1;M,\o)}
\newcommand{\modJplus}[1]{\mathcal{J}^+(#1;M,\o)}
\newcommand{\moduli}{\mathcal{M}}
\newcommand{\CP}[1]{\C P^{#1}}
\newcommand{\PSL}{PSL(2,\C)}
\newcommand{\Crit}{\mathrm{Crit}}
\newcommand{\im}{\mathrm{im}}
\newcommand{\partials}[1]{\frac{\partial #1}{\partial s}}
\newcommand{\partialt}[1]{\frac{\partial #1}{\partial t}}
\newcommand{\gen}[1]{\mathrm{<}#1\mathrm{>}}
\newcommand{\pss}{\cite{piunikhin-salamon-schwarz94}}

\newcommand{\ruled}{\mathbb{F}}
\newcommand{\grad}{\nabla}
\newcommand{\proj}{\mathbb{P}}
\newcommand{\pr}{\mathrm{pr}}
\newcommand{\mySp}{\mathrm{Sp}}
\newcommand{\Hom}{\mathrm{Hom}}
\newcommand{\eqclass}{S}
\newcommand{\tildeo}{\widetilde{\O}}

\standardspacing
\title[$\pi_1$ of symplectic automorphism groups]{$\pi_1$ of symplectic
automorphism groups and invertibles in quantum homology rings}
\date{Revised version, April 30, 1997.}
\thanks{Supported by a TMR grant from the European Community.}

\begin{document}
\maketitle
\newcommand{\hatj}{\hat{J}}
\newcommand{\jhat}{\hatj}
\section{Introduction}

The aim of this paper is to establish a connection between the
topology of the automorphism group of a symplectic manifold
$\mo$ and the quantum product on its homology. More precisely,
we assume that $M$ is closed and connected, and consider the group $\Ham\mo$
of Hamiltonian automorphisms with the $\smooth$-topology.
$\Ham\mo$ is a path-connected subgroup of the symplectic automorphism
group $\Aut\mo$; if $H^1(M,\R) = 0$, it is the connected component
of the identity in $\Aut\mo$. We introduce a homomorphism $q$
from a certain extension of the fundamental group $\pi_1(\Ham\mo)$ to the
group of invertibles in the quantum homology ring of $M$.
This invariant can be used to detect nontrivial elements
in $\pi_1(\Ham\mo)$. For example, consider 
$M = S^2 \times S^2$ with the family of product structures
$\o_\lambda = \lambda(\o_{S^2} \times 1) + 1 \times {\o_{S^2}}$, 
$\lambda > 0$. This example has been studied by Gromov \cite{gromov85},
McDuff \cite{mcduff87} and Abreu \cite{abreu97}. McDuff
showed that for $\lambda \neq 1$, $\pi_1(\Ham(M,\omega_\lambda))$ contains
an element of infinite order. This result can be recovered
by our methods. In a less direct way, the existence of $q$ 
imposes topological restrictions on all
elements of $\pi_1(\Ham\mo)$. An example
of this kind of reasoning can be found in section
\ref{sec:an-application}; a more important one will appear 
in forthcoming work by Lalonde, McDuff and Polterovich.

To define $q$, we will use the general relationship
between loops in a topological group and bundles over
$S^2$ with this group as structure group. In our case,
a smooth map $g: S^1 \longrightarrow \Aut\mo$ determines
a smooth fibre bundle $E_g$ over $S^2$ with a family
$\O_g = (\O_{g,z})_{z \in S^2}$ of symplectic structures
on its fibres. If $g(S^1)$ lies in $\Ham\mo$, the
$\O_{g,z}$ are restrictions of a closed $2$-form on 
$E_g$. We will call a pair $(E,\O)$ with this property a
{\em Hamiltonian fibre bundle}. 

Witten \cite{witten88} proposed to define an invariant of such a bundle
$(E,\O)$ over $S^2$ with fibre $\mo$ in the
following way: choose a positively oriented complex structure $j$ on $S^2$
and an almost complex structure $\hatj$ on $E$ (compatible
with the symplectic structure on each fibre) such that
the projection $\pi: E \longrightarrow S^2$ is $(\hatj,j)$-linear.
For generic $\hatj$, the space $\sec(j,\hatj)$ of pseudoholomorphic sections 
of $\pi$ is a smooth finite-dimensional manifold, with a natural evaluation
map $\ev_z: \sec(j,\hatj) \longrightarrow E_z$ for $z \in S^2$. 
After fixing a symplectic isomorphism 
$i : \mo \longrightarrow (E_z,\O_z)$ for some
$z$, $\ev_z$ should define a homology class in $M$. 
This rough description does not take into account the lack of
compactness of $\sec(j,\hatj)$. Part of this problem is due to
`bubbling': since we deal with it
using the method of Ruan and Tian in \cite{ruan-tian94}, 
we have to make the following

\begin{assumption}[$\mathrm{\mathbf{W^+}}$] 
$\mo$ satisfies one of the following
conditions:
\renewcommand{\theenumi}{(\alph{enumi})}
\renewcommand{\labelenumi}{\theenumi}
\begin{enumerate}
\item there is a $\lambda \geq 0$ such that $\o(A) = \lambda c_1(A)$
for all $A \in \hgroup$;
\item $c_1 | \hgroup = 0$;
\item the minimal Chern number $N \geq 0$ (defined by
$c_1(\hgroup) = N\Z$) is at least $n-1$, where $2n = \dim M$.
\end{enumerate}
\end{assumption}
Here $c_1$ stands for $c_1(TM,\o)$. An equivalent assumption is
\begin{equation} \label{eq:wplus-equivalent}
A \in \hgroup, 2-n \leq c_1(A) < 0 \Longrightarrow \o(A) \leq 0.
\end{equation}
Therefore {\wplus} is more restrictive than the notion of weak
monotonicity (see e.g. \cite{hofer-salamon95}) in which $2-n$ is
replaced by $3-n$. In particular, the Floer homology $HF_*\mo$ 
\cite{hofer-salamon95} and the quantum cup-product \cite{ruan-tian94} 
are well-defined for manifolds which satisfy {\wplus}. Note that
all symplectic four-manifolds belong to this class. Recent work on
Gromov-Witten invariants seems to indicate that the restriction {\wplus}
might be removed at the cost of introducing rational coefficients,
but this will not be pursued further here.

The other aspect of the non-compactness of $\sec(j,\hatj)$ is
that it may have infinitely many components corresponding
to different homotopy classes of sections of $E$. By separating
these components, one obtains infinitely many invariants
of $(E,\O)$. These can be arranged 
into a single element of the quantum homology group $\quantum$, which is
the homology of $M$ with coefficients in the Novikov ring
$\Lambda$. This element is normalized by the choice of a
section of $E$; it depends
on this section only up to a certain equivalence relation, which
we call {\em $\Gamma$-equivalence}. If $\eqclass$ is a 
$\Gamma$-equivalence class of sections, we denote the invariant obtained 
in this way by $Q(E,\O,S) \in \quantum$.

The bundles $(E_g,\O_g)$ for a Hamiltonian loop $g$ do not
come with a naturally preferred $\Gamma$-equivalence class of sections.
Therefore we introduce an additional piece of data:
$g$ acts on the free loop space $\Lambda M = \smooth(S^1,M)$ by
\begin{equation} \label{eq:action}
(g \cdot x)(t) = g_t(x(t)).
\end{equation}
Let $\loops$ be the connected component of $\Lambda M$ containing
the constant loops. There is an abelian covering $p: \tloops \longrightarrow
\loops$ such that $\Gamma$-equivalence classes of sections of $E_g$
correspond naturally to lifts of the action of $g$ to $\tloops$.
Let $\eqclass_{\tilde{g}}$ be the equivalence class corresponding to a
lift $\tilde{g}: \tloops \longrightarrow \tloops$. We define
\[
q\pair = Q(E_g,\O_g,\eqclass_{\tilde{g}}).
\] 
Let $G$ be the group of smooth loops $g: S^1 \longrightarrow \Ham\mo$
such that $g(0) = \id$, and $\group$ the group of pairs $\pair$. 
We will use the $\smooth$-topology on $G$, and a topology on $\group$
such that the homomorphism $\group \longrightarrow G$ which
`forgets' $\tilde{g}$ is continuous with discrete kernel.

The covering $\tloops$ was originally introduced by Hofer and
Salamon in their definition of $HF_*\mo$ \cite{hofer-salamon95}.
They set up Floer homology as a formal analogue of Novikov homology
for this covering. Now the $\group$-action on $\tloops$ commutes
with the covering transformations. 
Pursuing the analogy, one would expect an induced action of
$\group$ on Floer homology. If we assume that {\wplus} holds, 
this picture is correct: there are induced maps
\[
HF_*\pair: \floer \longrightarrow \floer
\]
for $(g,\tilde{g}) \in \tilde{G}$. These maps are closely related
to $q(g,\tilde{g})$. The relationship involves
the `pair-of-pants' product
\[
\pairofpants: \floer \times \floer \longrightarrow \floer
\]
and the canonical isomorphism $\Psi^+:
\quantum \longrightarrow \floer$ of Piunikhin, Salamon
and Schwarz \cite{piunikhin-salamon-schwarz94},\cite{schwarz97}.

\begin{bigtheorem} \label{bigth:comparison}
For any $\pair \in \group$ and $b \in \floer$, 
\[
HF_*\pair(b) = \Psi^+(q\pair) \pairofpants b. 
\] \end{bigtheorem}

This formula provides an alternative approach to
$q$ and $Q$, and we will use it to derive 
several properties of these invariants. Let $\qp$ be the
product on $\quantum$ obtained from the quantum cup-product 
by Poincar{\'e} duality. We will call $\qp$ the
quantum intersection product. The quantum homology ring
$(\quantum,\qp)$ is a ring with unit and commutative in the
usual graded sense.

\begin{bigcorollary} \label{bigth:invertible}
For any Hamiltonian fibre
bundle $(E,\O)$ over $S^2$ with fibre $\mo$ 
and any $\Gamma$-equivalence class $\eqclass$ of sections
of $E$, $Q(E,\O,\eqclass)$ is an invertible element of
$(\quantum,\qp)$. \end{bigcorollary}

By definition, $Q(E,\O,\eqclass)$ is homogeneous and even-dimensional,
that is, it lies in $QH_{2i}\mo$ for some $i \in \Z$. We will
denote the group (with respect to $\qp$) of 
homogeneous even-dimensional invertible elements of
$\quantum$ by $\inv$. 

\begin{bigcorollary} \label{bigth:homomorphism}
$q(g,\tilde{g})$ depends only
on $[g,\tilde{g}] \in \pi_0(\group)$, and
\[
q: \pi_0(\tilde{G}) \longrightarrow \inv
\]
is a group homomorphism. \end{bigcorollary}

Let $\Gamma$ be the group of covering transformations of
$\tloops$. Since the kernel of the homomorphism
$\group \longrightarrow G$ consists of the pairs
$(\id,\gamma)$ with $\gamma \in \Gamma$, there is an exact
sequence
\begin{equation} \label{eq:pizero}
\cdots \longrightarrow 
\Gamma \longrightarrow \pi_0(\tilde{G}) \longrightarrow
\pi_0(G) \longrightarrow 1.
\end{equation}
On the other hand, for every $\gamma \in \Gamma$ there is an element
$[M] \otimes \gen{\gamma} \in QH_*\mo$ which is easily seen to be invertible.
The map $\tau: \Gamma \longrightarrow \inv$ defined in this way is an
injective homomorphism.

\begin{bigproposition} \label{bigth:trivial}
For $\gamma \in \Gamma$, $q(\id,\gamma) = \tau(\gamma)$.
\end{bigproposition}

It follows that there is a unique homomorphism $\bar{q}$
such that the diagram
\begin{equation} \label{diag:basic}
\xymatrix{
\Gamma \ar[r] \ar@{=}[d] & \pi_0(\group) \ar[r] \ar[d]^{q}
& \pi_0(G) \ar[d]^{\bar{q}}\\
\Gamma \ar[r]^-{\tau} & \inv \ar[r] & \inv/\tau(\Gamma)
}
\end{equation}
commutes. $\bar{q}$ is more interesting for applications
to symplectic geometry than $q$ itself, since
$\pi_0(G) \iso \pi_1(\Ham\mo)$.

In the body of the paper, we proceed in a different order than
that described above. The next section contains the 
definition of $\group$ and of Hamiltonian fibre bundles.
Section \ref{sec:floer-homology} is a review (without proofs)
of Floer homology. Our definition is a variant of
that in \cite{hofer-salamon95}, the difference being that we
allow time-dependent almost complex structures. This makes
it necessary to replace weak monotonicity by the condition {\wplus}. 
The following three sections contain the definition and basic properties 
of the maps $HF_*\pair$. An argument similar to the use of
`homotopies of homotopies' in the definition of Floer homology
shows that $HF_*\pair$ depends only on $[g,\tilde{g}] \in
\pi_0(\group)$. Moreover, we prove that $HF_*\pair$ is an automorphism
of $HF_*\mo$ as a module over itself with the pair-of-pants product.
This second property is considerably simpler: after a slight modification
of the definition of $\pairofpants$, the corresponding equation holds
on the level on chain complexes. $Q(E,\O,\eqclass)$ is defined in
section \ref{sec:holomorphic-sections}, and section
\ref{sec:gluing} describes the gluing argument which
establishes the connection between this invariant and the
maps $HF_*\pair$. After that, we prove the results stated
above. The grading of $q(g,\tilde{g})$ is determined by a
simple invariant $I\pair$ derived from the
first Chern class. This relationship is exploited in section 
\ref{sec:an-application} to obtain some vanishing results for this
invariant. The final section contains two explicit computations of 
$q(g,\tilde{g})$, one of which is the case of $S^2 \times S^2$ 
mentioned at the beginning.

\subsubsection*{Acknowledgments} 
Discussions with S. Ag\-ni\-ho\-tri,
S. Do\-nald\-son and V. Mu\-noz were helpful in preparing
this paper. I have profited from the referee's comments
on an earlier version. I am particularly indebted to 
L. Polterovich for explaining to me the point of view 
taken in section \ref{sec:the-action} and for his encouragement. 
Part of this paper was written during a stay at the Universit{\'e}
de Paris-Sud (Orsay) and the Ecole Polytechnique.

\renewcommand{\theenumi}{(\roman{enumi})}
\renewcommand{\labelenumi}{\theenumi}
\section{\label{sec:elementary} Hamiltonian loops and fibre bundles}

Throughout this paper, $\mo$ is a closed connected symplectic
manifold of dimension $2n$. We will write $c_1$ for $c_1(TM,\o)$.
As in the Introduction, $\Lambda M = \smooth(S^1,M)$ denotes the
free loop space of $M$, $\loops \subset \Lambda M$ the
subspace of contractible loops, $\Ham\mo$ the group of Hamiltonian
automorphisms and $G$ the group of smooth
based loops in $\Ham\mo$. We use the $\smooth$-topology
on $\Ham\mo$ and $G$.

\begin{lemma} \label{th:smooth-loops} The canonical homomorphism
$\pi_0(G) \longrightarrow \pi_1(\Ham\mo)$
is an isomorphism. \end{lemma}

This follows from the fact that any continuous loop in 
$\Ham\mo$ can be approximated by a smooth loop. Note that it
is unknown whether $\Ham\mo$ with the $\smooth$-topology is
locally contractible for all $\mo$ (see the discussion on p. 321 of
\cite{mcduff-salamon}). However, there is a neighbourhood
$U \subset \Ham\mo$ of the identity such that each path component
of $U$ is contractible. This weaker property is sufficient to
prove Lemma \ref{th:smooth-loops}.

$G$ acts on $\Lambda M$ by \eqref{eq:action}.

\begin{lemma} \label{th:contractible-component}
$g(\loops) = \loops$ for every $g \in G$. \end{lemma}

The proof uses an idea which is also the
starting point for the definition of induced maps on Floer homology. 
Let $\H = \smooth(S^1 \times M, \R)$
be the space of periodic Hamiltonians\footnote{We 
identify $S^1 = \R/\Z$ throughout.}. The perturbed action
one-form $\alpha_H$ on $\Lambda M$ associated to $H \in \H$ is
\[
\alpha_H(x)\xi = \int_{S^1} \omega(\dot{x}(t) - X_H(t,x(t)),\xi(t)) dt
\]
where $X_H$ is the time-dependent Hamiltonian vector field
of $H$. The zero set of $\alpha_H$ consists of the $1$-periodic solutions of
\begin{equation} \label{eq:hamilton}
\dot{x}(t) = X_H(t,x(t)). 
\end{equation}
We say that a Hamiltonian $K_g \in \H$ generates $g \in G$ if
\begin{equation} \label{eq:generates}
\frac{\partial g_t}{\partial t}(y) = X_{K_g}(t,g_t(y)).
\end{equation}

\begin{lemma} \label{th:hg} Let $K_g$ be a Hamiltonian
which generates $g \in G$. For every $H \in \H$, define $H^g \in \H$ by
\[
H^g(t,y) = H(t,g_t(y)) - K_g(t,g_t(y)).
\]
Then $g^*\alpha_H = \alpha_{H^g}$. \end{lemma}

The proof is straightforward.

\proof[Proof of Lemma \ref{th:contractible-component}]
Assume that $g(\loops)$ is a connected component of $\Lambda M$
distinct from $\loops$. If $H$ is small, $\alpha_H$ has no
zeros in $g(\loops)$ and by Lemma \ref{th:hg}, $\alpha_{H^g}$
has no zeros in $\loops$. This contradicts the Arnol'd 
conjecture (now a theorem)
which guarantees the existence of at least one contractible
$1$-periodic solution of \eqref{eq:hamilton} for every 
Hamiltonian. \qed

Since we will only use this for manifolds satisfying
{\wplus}, we do not really need the Arnol'd conjecture in full
generality, only the versions established in \cite{hofer-salamon95} and
\cite{ono94}.

The space $\tloops$ introduced in \cite{hofer-salamon95} is
defined as follows: consider pairs
$(v,x) \in \smooth(D^2,M) \times \loops$ such that
$x = v|\partial D^2$. $\tloops$ is the set of equivalence classes
of such pairs with respect to the following relation:
$(v_0,x_0) \sim (v_1,x_1)$ if $x_0 = x_1$ and
$\o(v_0 \glue \overline{v_1}) = 0$, $c_1(v_0 \glue \overline{v_1}) = 0$.
Here $v_0 \glue \overline{v_1}: S^2 \longrightarrow M$ is the map
obtained by gluing together $v_0,v_1$ along the boundaries.
$p: \tloops \longrightarrow \loops$, $p(v,x) = x$, is a covering
projection for the obvious choice of topology on $\tloops$.
Define $\Gamma = \hgroup/\hgroup_0$, where $\hgroup_0$ is the
subgroup of classes $a$ such that $\o(a) = 0$, $c_1(a) = 0$.
$\Gamma$ can also be defined as a quotient of
$\im(\hgroup \longrightarrow H_2(M;\Z))$.
Therefore the choice of base point for $\pi_2(M)$ is irrelevant,
and any $A: S^2 \longrightarrow M$ determines a class $[A] \in \Gamma$.
Clearly $\o(A)$ and $c_1(A)$ depend only on $[A]$. By an abuse
of notation, we will write $\o(\gamma)$ and $c_1(\gamma)$ for 
$\gamma \in \Gamma$. $\Gamma$ is the covering group of $p$. 
It acts on $\tloops$ by `gluing in spheres': 
$[A] \cdot [v,x] = [A \glue v,x]$ 
(see \cite{hofer-salamon95} for a precise description).

\begin{lemma} \label{th:lift-action}
The action of any $g \in G$ on $\loops$ can
be lifted to a homeomorphism of $\tloops$. \end{lemma}

\proof Since $\tloops$ is a connected covering, it 
is sufficient to show that the action of $g$ on $\loops$
preserves the set of smooth maps $S^1 \longrightarrow \loops$ which can 
be lifted to $\tloops$. Such a map is given by a
$B \in \smooth(S^1 \times S^1, M)$ with $\o(B) = 0$, $c_1(B) = 0$.
Its image under $g$ is given by $B'(s,t) = g_t(B(s,t))$. Now
$\o(B') = \o(B)$ because $(B')^*\o = B^*\o + d\theta$, where
$\theta(s,t) = K_g(t,g_t(B(s,t))) dt$ for a Hamiltonian $K_g$
as in \eqref{eq:generates}. Similarly, $c_1(B') = c_1(B)$ because
there is an isomorphism $D: B^*TM \longrightarrow (B')^*TM$ of
symplectic vector bundles, given by $D(s,t) = Dg_t(B(s,t))$. \qed

\begin{definition} \label{def:group} $\group \subset
 G \times \mathrm{Homeo}(\tloops)$ is the subgroup of pairs
$\pair$ such that $\tilde{g}$ is a lift of the $g$-action
on $\loops$. \end{definition}

We give $\group$ the topology induced from the $\smooth$-topology
on $G$ and the topology of pointwise convergence on 
$\mathrm{Homeo}(\tloops)$. This makes $\group$ into a topological
group, essentially because a lift $\tilde{g}$ of a given $g$ is
determined by the image of a single point. The projection
$\group \longrightarrow G$ is onto by Lemma \ref{th:lift-action},
and since its kernel consists of the pairs $(\id,\gamma)$ with
$\gamma \in \Gamma$, there is an exact sequence
\[
1 \longrightarrow \Gamma \longrightarrow \group
\longrightarrow G \longrightarrow 1
\]
of topological groups, with $\Gamma$ discrete.

A point $c = [v,x] \in \tloops$ determines a preferred homotopy
class of trivializations of the symplectic vector bundle
$x^*(TM,\o)$. This homotopy class consists of the maps
$\tau_c: x^*TM \longrightarrow S^1 \times (\R^{2n},\o_0)$ which
can be extended over $v^*TM$, and it is independent of the choice
of the representative $(v,x)$ of $c$. For $\pair \in \group$,
\[
l(t) = \tau_{\tilde{g}(c)}(t) Dg_t(x(t)) \tau_c(t)^{-1} \qquad (t \in S^1)
\]
is a loop in $\mySp(2n,\R)$. Up to homotopy, it does not depend on $c$ and 
on the trivializations. We define the `Maslov index' $I\pair \in \Z$ by 
$I\pair = \mathrm{deg}(l)$, where 
$\mathrm{deg}: H_1(\mySp(2n,\R)) \longrightarrow \Z$ is the 
standard isomorphism induced by the determinant on
$U(n) \subset \mySp(2n,\R)$.

\begin{lemma} \label{th:easy-maslov} $I\pair$ depends only
on $[g,\tilde{g}] \in \pi_0(\group)$. The map 
$I: \pi_0(\group) \longrightarrow \Z$ is a homomorphism,
and $I(\id,\gamma) = c_1(\gamma)$ for all $\gamma \in \Gamma$.
\end{lemma}

We omit the proof. Because of \eqref{eq:pizero}, it follows
that $I(g,\tilde{g}) \text{ mod } N$ depends only on $g$.
In this way, we recover a familiar (see e.g. \cite[p. 80]{weinstein89}) 
invariant
\[
\bar{I}: \pi_0(G) = \pi_1(\Ham\mo) \longrightarrow \Z/N\Z.
\]

Now we turn to symplectic fibre bundles. A smooth fibre
bundle $\pi: E \longrightarrow B$ together with a smooth family
$\O = (\O_b)_{b \in B}$ of symplectic forms on its fibres
is called a symplectic fibration;
a symplectic fibre bundle is a symplectic fibration which is
locally trivial. Note that $\O$ defines a symplectic structure
on the vector bundle $TE^v = \ker(D\pi) \subset TE$.

Here we consider only the case $B = S^2$.
It is convenient to think of $S^2$ as $D^+ \cup_{S^1} D^-$,
where $D^+$, $D^-$ are closed discs. Fix a point $z_0 \in D^-$.
A symplectic fibre bundle $(E,\O)$ over $S^2$ with a fixed isomorphism
$i: \mo \longrightarrow (E_{z_0},\O_{z_0})$ will be called
a symplectic fibre bundle {\em with fibre $\mo$}; we will
frequently omit $i$ from the notation.

Let $g$ be a smooth loop in $\Aut\mo$. The `clutching' construction
produces a symplectic fibre bundle $(E_g,\O_g)$ over
$S^2$ by gluing together the trivial fibre
bundles $D^{\pm} \times (M,\omega)$ using
\[ \begin{gathered}
\phi_g: \partial D^+ \times M \longrightarrow
\partial D^- \times M,\\
\phi_g(t,y) = (t,g_t(y)).
\end{gathered} \]
In an obvious way, $(E_g,\O_g)$ is a symplectic fibre bundle
with fibre $\mo$. 

\begin{convention} \label{th:orientation-convention}
We have identified $\partial D^+,\partial D^-$
with $S^1$ and used these identifications to glue $D^+$ and
$D^-$ together. We orient $D^+,D^-$ in such a way
that the map $S^1 \longrightarrow \partial 
D^+$ preserves orientation while the one $S^1 \longrightarrow
\partial D^-$ reverses it; this induces an orientation
of $S^2$.
\end{convention}

We will use this construction only for loops with $g(0) = \id$. 
By a standard argument, it provides
a bijection between elements of $\pi_1(\Aut\mo)$ and isomorphism classes
of symplectic fibre bundles over $S^2$ with fibre $\mo$. 

\begin{defn} A symplectic fibre bundle $(E,\O)$ over a
surface $B$ is called {\em Hamiltonian} if there
is a closed two-form $\tildeo$ on $E$ such that
$\tildeo|E_b = \O_b$ for all $b \in B$. \end{defn}

For $B = S^2$, these are precisely
the bundles corresponding to Hamiltonian loops:

\begin{proposition} \label{th:hamiltonian-bundles}
$(E_g,\O_g)$ is Hamiltonian iff $[g]$ lies in the subgroup\\
$\pi_1(\Ham\mo) \subset \pi_1(\Aut\mo)$. \end{proposition}

\proof Let $g \in G$ be a loop generated by $K_g \in \H$.
We denote the pullback of $\o$ to $D^\pm \times M$ by $\o^{\pm}$.
Choose a $1$-form $\theta$ on $D^+ \times M$ such that
$\theta(t,y) = K_g(t,g_t(y)) dt$ for $t \in S^1$ and
$\theta | \{z\} \times M = 0$ for all $z \in D^+$. There is
a continuous $2$-form
$\tildeo$ on $E_g$ such that $\tildeo | D^+ \times M =
\o^+ + d\theta$ and $\tildeo | D^- \times M = \o^-$, 
and for a suitable choice of $\theta$ it is smooth. 
$\tildeo$ is closed and extends the family
$\O_g$, which proves that $(E_g,\O_g)$ is Hamiltonian.

Conversely, let $g$ be a based loop in $\Aut\mo$ such that
$(E_g,\O_g)$ is Hamiltonian with $\tildeo \in \Omega^2(E_g)$. 
We will use the exact sequence
\begin{equation} \label{eq:banyaga}
1 \longrightarrow \pi_1(\Ham\mo) \longrightarrow
\pi_1(\Aut\mo) \stackrel{F}{\longrightarrow} H^1(M,\R)
\end{equation}
where $F$ is the flux homomorphism (see \cite{banyaga78} or
\cite[Corollary 10.18]{mcduff-salamon96}). Choose a point
$z_+ \in D^+$. For $\lambda \in \smooth(S^1,M)$, consider the maps
\begin{align*}
T_{\lambda}: S^1 \times S^1 &\longrightarrow M, \quad 
T_{\lambda}(r,t) = g_t(\lambda(r)),\\
T_{\lambda}^-: S^1 \times S^1 &\longrightarrow D^- \times M, \quad 
T_{\lambda}^-(r,t) = (z_0,g_t(\lambda(r))) \text{ and } \\
T_{\lambda}^+: S^1 \times S^1 &\longrightarrow D^+ \times M, \quad 
T_{\lambda}^+(r,t) = (z_+,\lambda(r)).
\end{align*}
By definition
$\left<F(g),[\lambda]\right> = \o(T_\lambda) = \tildeo(T_{\lambda}^-)$.
Since $T_{\lambda}^-$ and $T_{\lambda}^+$ are homotopic in $E_g$,
$\tildeo(T_{\lambda}^-) = \tildeo(T_{\lambda}^+)$; but 
clearly $\tildeo(T_{\lambda}^+) = 0$. Therefore $F(g) = 0$, and
by \eqref{eq:banyaga}, $[g] \in \pi_1(\Ham\mo)$. \qed

Let $(E,\O)$ be a Hamiltonian fibre bundle over $S^2$, 
with $\tildeo \in \Omega^2(E)$ as above.
We say that two continuous sections $s_0$, $s_1$ of $E$ are
{\em $\Gamma$-equivalent} if $\tildeo(s_0) = \tildeo(s_1)$ and
$c_1(TE^v,\O)(s_0) = c_1(TE^v,\O)(s_1)$. Using the exact sequence
\begin{equation} \label{eq:homotopy-sequence}
\cdots \longrightarrow \pi_2(E_{z_0}) 
\longrightarrow \pi_2(E) \longrightarrow
\pi_2(S^2) \longrightarrow \cdots,
\end{equation}
it is easy to see that this equivalence relation is independent of the
choice of $\tildeo$. Let $S$ be the $\Gamma$-equivalence class of
a section $s$. By an abuse of notation, we write $\tildeo(S) =
\tildeo(s)$ and $c_1(TE^v,\O)(S) = c_1(TE^v,\O)(s)$.

\begin{lemma} \label{th:equivalence-classes}
Let $(E,\O)$ be a Hamiltonian fibre bundle over $S^2$ with fibre $\mo$
and $\tildeo \in \Omega^2(E)$ a closed extension of $\O$.
\begin{enumerate} \item \label{item:cont-sec}
$E$ admits a continuous section.
\item \label{item:gamma-one}
For two $\Gamma$-equivalence classes $S_0$, $S_1$ of sections of $E$,
there is a unique $\gamma \in \Gamma$ such that
\begin{equation} \label{eq:oc-difference}
\begin{aligned}
\tildeo(S_1) &= \tildeo(S_0) + \o(\gamma),\\
c_1(TE^v,\O)(S_1) &= c_1(TE^v,\O)(S_0) + c_1(\gamma).
\end{aligned}
\end{equation}
$\gamma$ is independent of the choice of $\tildeo$.
\item \label{item:gamma-two}
Conversely, given a $\Gamma$-equivalence class $S_0$ and
$\gamma \in \Gamma$, there is a unique $\Gamma$-equivalence class
$S_1$ such that \eqref{eq:oc-difference} holds. 
\end{enumerate} \end{lemma}

\proof \ref{item:cont-sec} $(E,\O)$ is isomorphic to $(E_g,\O_g)$
for some $g: S^1 \longrightarrow \Aut\mo$. Since it is Hamiltonian,
we can assume that $g \in G$. Choose a point $y \in M$; by Lemma 
\ref{th:contractible-component}, there
is a $v \in \smooth(D^-,M)$ such that $v(t) = g_t(y)$ for
$t \in S^1 = \partial D^-$. To obtain a section of $E_g$, glue together
$s^+: D^+ \longrightarrow D^+ \times M$, $s^+(z) = (z,y)$ and
$s^-: D^- \longrightarrow D^- \times M$, $s^-(z) = (z,v(z))$. \\
\ref{item:gamma-one} is an easy consequence of 
\eqref{eq:homotopy-sequence}. \\
\ref{item:gamma-two} Let $i: \mo \longrightarrow (E_{z_0},\O_{z_0})$ be the
preferred symplectic isomorphism. Take $s_0 \in S_0$ and 
$A: S^2 \longrightarrow M$ such that $[A] = \gamma$ and 
$i(A(z)) = s_0(z_0)$ for some $z \in S^2$. A section $s_1$ with the 
desired properties can be obtained by gluing together $s_0$ and $i(A)$. 
The uniqueness of $S_1$ is clear from the definition. \qed

In the situation of Lemma \ref{th:equivalence-classes}, we will
write $\gamma = S_1 - S_0$ and $S_1 = \gamma + S_0$.

\begin{defn} Let $(E,\O)$ be a Hamiltonian fibre bundle over $S^2$
and $S$ a $\Gamma$-equivalence class of sections of
$E$. We will call $(E,\O,S)$ a {\em normalized Hamiltonian fibre bundle}. 
\end{defn}

This notion is closely related to the group $\group$:
for $\pair \in \group$, choose a point $c \in \tloops$
and representatives $(v,x)$ of $c$ and $(v',x')$ of $\tilde{g}(c)$.
The maps
$s^+_{\tilde{g}}: D^+ \longrightarrow D^+ \times M$, 
$s^+_{\tilde{g}}(z) = (z,v(z))$ and
$s^-_{\tilde{g}}: D^- \longrightarrow D^- \times M$, 
$s^-_{\tilde{g}}(z) = (z,v'(z))$
define a section $s_{\tilde{g}}$ of $E_g$. To be precise, 
$s^+_{\tilde{g}}$ is obtained from $v$ by identifying $D^+$ with the
standard disc $D^2$ such that the orientation (see
\ref{th:orientation-convention}) is
preserved; in the case of $D^-$, one uses an orientation-reversing 
diffeomorphism. By comparing this with the definition of $I\pair$, 
one obtains
\begin{equation} \label{eq:maslov-and-chern}
c_1(TE_g^v,\O_g)(s_{\tilde{g}}) = -I(g,\tilde{g}).
\end{equation} 

\begin{lemma} \label{th:equivalence-class}
The $\Gamma$-equivalence class of $s_{\tilde{g}}$ is independent of the
choice of $c$ and of $v,v'$. \end{lemma}

We omit the proof. It follows that $\pair$ determines a normalized 
Hamiltonian fibre bundle $(E_g,\O_g,S_{\tilde{g}})$, where $S_{\tilde{g}}$
is the $\Gamma$-equivalence class of $s_{\tilde{g}}$. 
This defines a one-to-one correspondence between $\pi_0(\group)$ 
and isomorphism classes of normalized Hamiltonian fibre bundles 
with fibre $\mo$. We will only use the easier half of this correspondence:

\begin{lemma} \label{th:completeness}
Every normalized Hamiltonian fibre bundle with fibre $\mo$ is isomorphic
to $(E_g,\O_g,\eqclass_{\tilde{g}})$ for some 
$\pair \in \group$. \qed \end{lemma}

\newcommand{\floermetric}{\smooth_{\epsilon}}
\section{\label{sec:floer-homology} Floer homology}

From now on, we always assume that $\mo$ satisfies {\wplus}.
In \cite{hofer-salamon95}, weak monotonicity is used
to show that a generic $\o$-compatible almost complex structure
admits no pseudo-holomorphic spheres with negative Chern
number. Under the more restrictive condition
{\wplus}, this non-existence result can be extended
to families of almost complex structures depending on
$\leq\!3$ parameters. To make this precise, we need to introduce
some notation. Let $\J = (J_b)_{b \in B}$ be a smooth family 
of almost complex structures on $M$ parametrized by a manifold
$B$. $\J$ is called $\o$-compatible if every $J_b$ is
$\o$-compatible in the usual sense. We denote the space of
$\o$-compatible families by $\modJ{B}$.

The parametrized moduli space 
$\moduli^s(\J) \subset B \times \smooth(\CP{1},M)$ associated to
$\J \in \modJ{B}$ is the space of pairs 
$(b,w)$ such that $w$ is $J_b$-holomorphic and simple (not multiply 
covered). In the case of a single almost complex structure $J$, 
a $J$-holomorphic sphere $w$ is called regular if the linearization
of the equation $\bar{\partial}_J(u) = 0$ at $w$, which is
given by an operator
\[
D_J(w): \smooth(w^*TM) \longrightarrow \Omega^{0,1}(w^*(TM,J)),
\]
is onto (see \cite[Chapter 3]{mcduff-salamon}). Similarly,
$(b,w) \in \moduli^s(\J)$ is called regular if the extended operator
\begin{gather*}
\hat{D}_{\J}(b,w): T_bB \times
\smooth(w^*TM) \longrightarrow \Omega^{0,1}(w^*(TM,J_b)), \\
\hat{D}_{\J}(b,w)(\beta,W) =
D_{J_b}(w)W + \frac{1}{2} \left( D\J(b)\beta \right) \circ dw \circ i
\end{gather*}
is onto. Here $i$ is the complex structure on $\CP{1}$
and $D\J(b): T_bB \longrightarrow \smooth(\mathrm{End}(TM))$
is the derivative of the family $(J_b)$ at $b$.
$\J$ itself is called regular if all $(b,w) \in \moduli^s(\J)$ are regular,
and the set of regular families is denoted by
$\modJreg{B} \subset \modJ{B}$.

$\hat{D}_{\J}(b,w)$ is a Fredholm operator of index
$2n + 2c_1(w) + \dim B$. For $k \in \Z$, let $\moduli^s_k(\J)
\subset \moduli^s(\J)$ be the subspace of pairs $(b,w)$ with
$c_1(w) = k$. By applying the implicit function theorem, one
obtains

\begin{lemma} \label{lemma:implicit-function}
If $\J \in \modJreg{B}$, $\moduli^s_k(\J)$ is a smooth manifold of dimension
$2n + 2k + \dim B$ for all $k$. \end{lemma}

Assume that $B$ is compact and choose $\J_0 \in
\modJ{B}$. Let $U_\delta(\J_0) \subset \modJ{B}$ be a
$\delta$-ball around $\J_0$ with respect to Floer's
$\floermetric$-norm (see \cite{hofer-salamon95}
or \cite[p. 101--103]{schwarz95}).

\begin{theorem} \label{th:para-trans-theorem}
For sufficiently small $\delta>0$,
$\modJreg{B} \cap U_\delta(\J_0) \subset U_\delta(\J_0)$
has second category; in particular, it is $\floermetric$-dense.
\end{theorem}

This is a well-known result (see e.g. \cite[Theorem 3.1.3]{mcduff-salamon}).
It is useful to compare its proof with that of the basic transversality
theorem for pseudoholomorphic curves, which is the special case
$B = \{pt\}$: in that case, one shows first that the `universal moduli space'
\[
\moduli^{\mathrm{univ}} \stackrel{\pi}{\longrightarrow} \modJ{pt}
\]
is smooth. $\modJreg{pt}$ is the set of regular values of $\pi$,
which is shown to be dense by applying the Sard-Smale theorem
(we omit the details which arise from the use of $\floermetric$-spaces).
In the general case, a family $(J_b)_{b \in B}$ is regular iff the
corresponding map $B \longrightarrow \modJ{pt}$ is transverse to $\pi$.
Therefore the part of the proof which uses specific properties of
pseudoholomorphic curves (the smoothness of $\moduli^{\mathrm{univ}}$)
remains the same, but a different general result
has to be used to show that generic maps $B \longrightarrow \modJ{pt}$
are transverse to $\pi$.

A family $\J$ is called semi-positive if $\moduli^s_k(\J) = \emptyset$
for all $k < 0$; the semi-positive families form a subset
$\modJplus{B} \subset \modJ{B}$.

\begin{lemma} If $\dim B \leq 3$, $\modJreg{B} \subset \modJplus{B}$. 
\end{lemma}

\proof $\PSL$ acts freely on $\moduli_k^s(\J)$ for all $k$. 
If $\J \in \modJreg{B}$, the quotient is a smooth manifold and
\[
\dim \moduli_k^s(\J)/\PSL \leq 2n + 2k - 3.
\]
Therefore $\moduli_k^s(\J) = \emptyset$ for $k \leq 1-n$.
But for $2-n \leq k < 0$, $\moduli_k^s(\J) = \emptyset$
by \eqref{eq:wplus-equivalent}. \qed

\begin{corollary} \label{cor:no-negative-spheres}
If $B$ is compact and $\dim B \leq 3$, $\modJplus{B}$
is $\smooth$-dense in $\modJ{B}$. \end{corollary}

Later, we will also use a `relative'
version of Corollary \ref{cor:no-negative-spheres}, in which $B$ is 
non-compact and one considers families $\J$ with a fixed behaviour 
outside a relatively compact open subset of $B$. We omit the 
precise statement.

In contrast with the case of negative Chern number, holomorphic
spheres with Chern number $0$ or $1$ can occur in a family
$\J = (J_b)_{b \in B}$ even if $\dim B$ is small. When defining
Floer homology, special attention must be paid to them. For any
$k \geq 0$, let
$V_k(\J) \subset B \times M$ be the set of pairs $(b,y)$ such
that $y \in \im(w)$ for a non-constant $J_b$-holomorphic sphere $w$ with
$c_1(w) \leq k$. $V_k(\J)$ is the 
union of the images of the evaluation maps
\[
\moduli_j^s(\J) \times_{\PSL} \CP{1} \longrightarrow B \times M
\]
for $j \leq k$. For regular $\J$ and $j \leq 0$, 
the dimension of $\moduli_j^s(\J) \times_{\PSL} \CP{1}$ is
$\leq 2n + \dim B - 4$. Therefore $V_0(\J)$
is (loosely speaking) a codimension-$4$ subset of
$B \times M$. Similarly, $V_1(\J)$ has codimension $2$.

For any $H \in \H$, the pullback of $\alpha_H$ to $\tloops$
is exact: $p^*\alpha_H = da_H$ with
\[
a_H(v,x) = -\int_{D^2} v^*\o + \int_{S^1} H(t,x(t)) dt.
\]
Let $Z(\alpha_H) \subset \loops$ be the zero set of
$\alpha_H$ and $\Crit(a_H) = p^{-1}(Z(\alpha_H))$ the
set of critical points of $a_H$. $c = [v,x] \in \Crit(a_H)$
is nondegenerate iff $x$ is a regular zero of $\alpha_H$,
that is, if the linearization of \eqref{eq:hamilton} at
$x$ has no nontrivial solutions. We denote the Conley-Zehnder
index of such a critical point by $\mu_H(c) \in \Z$.

\begin{convention} Let $H$ be a time-independent Hamiltonian,
$y \in M$ a nondegenerate critical point of $H$ and $c = [v,x] \in \tloops$ 
the critical point of $a_H$ represented by the constant maps $v,x \equiv y$.
Our convention for $\mu_H$ is that for small $H$, $\mu_H(c)$ equals 
the Morse index of $y$. This differs from the convention in
\cite{salamon-zehnder92} by a constant $n$. \end{convention}

Assume that all critical points of $a_H$ are nondegenerate, and
let $\Crit_k(a_H)$ be the set of critical points
of index $k$. The Floer chain group $CF_k(H)$ is the group of formal sums
\[
\sum_{c \in \Crit_k(a_H)} m_c \gen{c}
\]
with coefficients $m_c \in \Z/2$ such that
$
\{ c \in \Crit_k(a_H) \, | \, m_c \neq 0, a_H(c) \geq C \}
$
is finite for all $C \in \R$.

A family $\J = (J_t)_{t \in S^1} \in \modJ{S^1}$ defines a
Riemannian metric
\[
(\xi,\eta)_{\J} = \int_{S^1} \o(\xi(t),J_t \eta(t)) dt  
\]
on $\loops$. Let $\nabla_{\J} a_H$ be the gradient of $a_H$ with respect
to the pullback of $(\cdot,\cdot)_{\J}$ to $\tloops$. A smooth path 
$\tilde{u}: \R \longrightarrow \tloops$ is a flow line of
$-\nabla_{\J} a_H$ iff its projection to $\loops$
is given by a map $u \in \smooth(\R \times S^1,M)$ such that
\begin{equation} \label{eq:floer}
\partials{u} + J_t(u(s,t))\left(
\partialt{u} - X_H(t,u(s,t))\right) = 0
\quad \text{ for all } (s,t) \in \R \times S^1.
\end{equation}
If $u$ is a solution of \eqref{eq:floer} whose 
energy $E(u) = \int {\left|\partials{u}\right|}^2$ is
finite, there are $c_-,c_+ \in \Crit(a_H)$ such that
\[
\lim_{s \rightarrow \pm \infty} \tilde{u}(s) = c_\pm.
\]
Moreover, if the limits are nondegenerate, the linearization 
of \eqref{eq:floer} at $u$ is a differential operator
$D_{H,\J}(u): \smooth(u^*TM) \longrightarrow \smooth(u^*TM)$
whose Sobolev completion
\[
D_{H,\J}(u): W^{1,p}(u^*TM) \longrightarrow L^p(u^*TM)
\]
$(p>2)$ is a Fredholm operator with
index $\ind(u) = \mu_H(c_-) - \mu_H(c_+)$. 

For $c_-,c_+ \in \Crit(a_H)$, let $\moduli(c_-,c_+;H,\J)$
be the space of solutions of \eqref{eq:floer} which have
a lift $\tilde{u}: \R \longrightarrow \tloops$ with limits
$c_-,c_+$. The Floer differential
\[
\partial_k(H,\J): CF_k(H) \longrightarrow CF_{k-1}(H)
\]
is defined by the formula
\begin{equation} \label{eq:differential}
\partial_k(H,\J)(\gen{c_-}) = \sum_{c_+ \in \Crit_{k-1}(a_H)}
\#(\moduli(c_-,c_+;H,\J)/\R) \gen{c_+},
\end{equation}
extended in the obvious way to infinite linear combinations
of the $\gen{c_-}$ (from now on, this is to be understood 
for all similar formulae). $\R$ acts on $\moduli(c_-,c_+;H,\J)$
by $(s_0 \cdot u)(s,t) = u(s-s_0,t)$, and
$\#$ denotes counting the number
of points in a set mod $2$. In order for \eqref{eq:differential}
to make sense and define the `right' homomorphism, $(H,\J)$
has to satisfy certain conditions.

\begin{defn} \label{def:regular-pair}
$(H,\J) \in \H \times \modJreg{S^1}$ is called a
{\em regular pair} if
\begin{enumerate}
\item \label{item:one} every $[v,x] \in \Crit(a_H)$ is nondegenerate and
satisfies $(t,x(t)) \notin V_1(\J)$ for all $t \in S^1$;
\item \label{item:regul}
every solution $u: \R \longrightarrow \loops$ 
of \eqref{eq:floer} with finite energy
is regular (that is, the linearization of \eqref{eq:floer} at $u$ is onto);
\item \label{item:three}
if in addition $\ind(u) \leq 2$, $(t,u(s,t)) \notin V_0(\J)$
for all $(s,t) \in \R \times S^1$.
\end{enumerate} \end{defn}

These are essentially the conditions defining
$\mathcal{H}_{\mathrm{reg}}(J)$ in \cite{hofer-salamon95},
adapted to the case of time-dependent almost complex
structures. \ref{item:regul} implies that $\moduli(c_-,c_+;H,\J)$
is a manifold of dimension $\mu_H(c_-) - \mu_H(c_+)$ for all
$c_-,c_+$. The proof of the main compactness result
\cite[Theorem 3.3]{hofer-salamon95} carries over to our
situation using Corollary \ref{cor:no-negative-spheres}. 
It follows that if $(H,\J)$ is a regular pair, the r.h.s. of
\eqref{eq:differential} is meaningful and defines differentials
$\partial_k(H,\J)$ such that $\partial_{k-1}(H,\J) \partial_k(H,\J) = 0$.

\begin{definition} The Floer homology $HF_*(H,\J)$ of a regular
pair $(H,\J)$ is the homology of $(CF_*(H),\partial_*(H,\J))$. \end{definition}

The existence of regular pairs is ensured by the following Theorem,
which is a variant of \cite[Theorems 3.1 and 3.2]{hofer-salamon95}.

\begin{thm} \label{th:transv-floer} 
The set of regular pairs is $\smooth$-dense
in $\H \times \modJ{S^1}$. \end{thm}

So far, we have only used families of almost complex structures 
parametrized by $S^1$. For such families Corollary 
\ref{cor:no-negative-spheres} holds whenever $\mo$ is weakly
monotone; therefore, the stronger assumption {\wplus} has
not been necessary up to now. Two-parameter
families of almost complex structures occur first in the next step,
the definition of `continuation maps'.

A homotopy between regular pairs $(H^-,\J^-)$, $(H^+,\J^+)$ consists
of an $H \in \smooth(\R \times S^1 \times M,\R)$ and a
$\J \in \modJreg{\R \times S^1}$, such that $(H(s,t,\cdot), J_{s,t})$
is equal to $(H^-(t,\cdot),J^-_t)$ for $s \leq -1$ and to
$(H^+(t,\cdot),J^+_t)$ for $s \geq 1$. For $c_- \in \Crit(a_{H^-})$
and $c_+ \in \Crit(a_{H^+})$, let $\moduli^{\Phi}(c_-,c_+;H,\J)$ be the
space of solutions $u \in \smooth(\R \times S^1, M)$ of
\begin{equation} \label{eq:continuation}
\partials{u} + J_{s,t}(u)\left(\partialt{u} - X_H(s,t,u)\right) = 0
\end{equation}
which can be lifted to paths $\tilde{u}: \R \longrightarrow \tloops$
with limits $c_-,c_+$. A homotopy $(H,\J)$ is regular if every
solution of \eqref{eq:continuation} is regular (the linearization
is onto) and if for those solutions with index $\leq 1$,
\[
(s,t,u(s,t)) \notin V_0(\J) \text{ for all } (s,t) \in \R \times S^1.
\]
For a regular homotopy, the `continuation homomorphisms'
\[
\Phi_k(H,\J): CF_k(H^-,\J^-) \longrightarrow CF_k(H^+,\J^+)
\]
are defined by
\[
\Phi_k(H,\J)(\gen{c_-}) = \sum_{c_+ \in \Crit_k(a_{H^+})}
 \#\moduli^{\Phi}(c_-,c_+;H,\J) \gen{c_+}.
\]
This is a homomorphism of chain complexes and induces an isomorphism of
Floer homology groups, which we denote equally by $\Phi_*(H,\J)$.
For fixed $(H^{\pm},\J^{\pm})$, these isomorphisms are
independent of the choice of $(H,\J)$. Moreover, they
are functorial with respect to composition of homotopies.
Therefore we can speak of a well-defined Floer homology
$HF_*\mo$ independent of the choice of a regular pair
(see \cite{salamon-zehnder92} for a detailed discussion,
which carries over to our case with minor modifications).
The proof that the `continuation maps' are independent
of $(H,\J)$ involves three-parameter families of almost
complex structures, and hence uses the full strength of
Corollary \ref{cor:no-negative-spheres}.

\begin{remark} The definition used here leads to a Floer homology
which is canonically isomorphic to that defined in \cite{hofer-salamon95}.
The only difference is that we have allowed a wider range
of perturbations. The reason for this will become clear in the next
section. \end{remark}

Since $a_H = p^*\alpha_H$, $\Crit(a_H)$ is $\Gamma$-invariant.
This induces an action of $\Gamma$ on $CF_*(H)$ for every regular
pair $(H,\J)$. The action does not preserve the grading:
$\gamma$ maps $CF_k(H)$ to $CF_{k - 2c_1(\gamma)}(H)$, because
\begin{equation} \label{eq:grading-shift}
\mu_H(A \glue v,x) = \mu_H(v,x) - 2c_1(A)
\end{equation}
for $[v,x] \in \Crit(a_H)$ and $A: S^2 \longrightarrow M$
(see \cite{hofer-salamon95} or \cite[Proposition 5]{floer-hofer93};
the difference in sign to \cite{hofer-salamon95} is due to the
fact that they consider Floer cohomology rather than homology).
To recover a proper grading, one should assign to $\gamma \in \Gamma$ the
`dimension' $-2c_1(\gamma)$. We denote the subset of elements
of `dimension' $k$ by $\Gamma_k \subset \Gamma$.

\begin{defn} \label{def:novikov-ring}
For $k \in \Z$, let $\novi_k$ be the group of formal sums
\begin{equation} \label{eq:element-of-novikov-ring}
\sum_{\gamma \in \Gamma_k} m_\gamma \gen{\gamma}
\end{equation}
with $m_\gamma \in \Z/2$, such that $\{\gamma \in \Gamma_k \,
| \, m_\gamma \neq 0, \o(\gamma) \leq C\}$ is finite for all
$C \in \R$. The multiplication
\[
\gen{\gamma} \cdot \gen{\gamma'} = \gen{\gamma + \gamma'}
\]
can be extended to the infinite linear combinations
\eqref{eq:element-of-novikov-ring}, and this makes
$\novi = \bigoplus_k \novi_k$ into a commutative graded ring
called the Novikov ring of $\mo$. \end{defn}

For every $\gamma \in \Gamma$ and $c \in \tloops$, 
$a_H(\gamma \cdot c) = a_H(c) - \o(\gamma)$.
Therefore the $\Gamma$-action on $CF_*(H)$ extends naturally to a graded
$\novi$-module structure. Since the differentials $\partial_*(H,\J)$ and
the homomorphisms $\Phi_*(H,\J)$ are $\novi$-linear,
this induces a $\novi$-module structure on $HF_*(H,\J)$ and
on $HF_*\mo$. 

The {\em quantum homology} $\quantum$ is the graded $\novi$-module
defined by
\[
QH_k\mo = \bigoplus_{i + j = k} H_i(M;\Z/2) \otimes
\novi_j.
\]

\begin{theorem} \label{th:what-is-floer-homology}
$\floer$ is isomorphic to $\quantum$ as a
graded $\novi$-module. \end{theorem}

This was proved by Piunikhin, Salamon and Schwarz {\pss} (less general
versions had been obtained before by Floer \cite{floer88} and
Hofer-Salamon \cite{hofer-salamon95}); a detailed account of the
proof is in preparation \cite{schwarz97}. Their result is in fact
much stronger: it provides a canonical isomorphism
$\Psi^+: \quantum \longrightarrow \floer$ which also relates
different product structures. We will use this construction 
in section \ref{sec:gluing}.

\section{\label{sec:the-action} The $\group$-action on Floer homology}

We begin by considering the action of $G$ on $\modJ{S^1}$. For 
$\J = (J_t)_{t \in S^1} \in \modJ{S^1}$ and $g \in G$, define $\J^g =
(J_t^g)_{t \in S^1}$ by $J_t^g = Dg_t^{-1} J_t Dg_t$.

\begin{lemma} \label{th:naturality-of-regularity}
If $\J \in \modJreg{S^1}$, $\J^g \in \modJreg{S^1}$ for
all $g \in G$. \end{lemma}

\proof $w \in \smooth(\CP{1},M)$ is $J_t^g$-holomorphic iff
$w' = g_t(w)$ is $J_t$-holomorphic. Since
\begin{align*}
\hat{D}_{\J}(t,w')(1,\frac{\partial g_t}{\partial t}(w)) &=
\frac{\partial}{\partial t}
\left(\frac{1}{2} Dg_t \circ dw + \frac{1}{2} J_t \circ Dg_t \circ dw
\circ i\right)\\
&= \frac{1}{2} Dg_t \circ \frac{\partial J_t^g}{\partial t} \circ dw
\circ i,
\end{align*}
there is a commutative diagram
\begin{equation*} \xymatrix{
\R \times \smooth(w^*TM) \ar[rr]^{\hat{D}_{\J^g}(t,w)} \ar[d]^{\delta}
&& \Omega^{0,1}(w^*(TM,J_t^g)) \ar[d]^{\delta'}\\
\R \times \smooth((w')^*TM) \ar[rr]^{\hat{D}_{\J}(t,w')} && 
\Omega^{0,1}((w')^*(TM,J_t)).
} 
\end{equation*}
Here $\delta(\tau,W) = (\tau, Dg_t(w)W + \tau \frac{\partial
g_t}{\partial t}(w))$, and $\delta'$ maps a homomorphism
$\sigma: T\CP{1} \longrightarrow w^*TM$ to $\delta'(\sigma) =
Dg_t(w) \circ \sigma$. If $w$ is simple, $\hat{D}_\J(t,w')$
is onto by assumption, and since $\delta,\delta'$ are isomorphisms,
$\hat{D}_{\J^g}(t,w)$ is onto as well. \qed

The same argument shows that $\hat{D}_{\J^g}(t,w)$ and
$\hat{D}_\J(t,w')$ have the same index; therefore
\begin{equation} \label{eq:g:v}
V_k(\J^g) = \{(t,y) \in S^1 \times M \; | \; (t,g_t(y)) \in V_k(\J)\}
\end{equation}
for all $k \geq 0$. The metric on $\loops$ associated to $\J^g$ is
\begin{equation} \label{eq:g:metric}
(\cdot,\cdot)_{\J^g} = g^*(\cdot,\cdot)_\J.
\end{equation}

\begin{remark} The $G$-action does not preserve the subspace
of families $(J_t)$ such that $J_t = J_0$ for all $t$.
This is the reason why we have allowed $t$-dependent almost complex
structures in the definition of Floer homology. \end{remark}

For $(H,\J) \in \H \times \modJ{S^1}$ and $g \in G$, we call
the pair $(H^g,\J^g)$, where $H^g$ is as in Lemma \ref{th:hg},
the {\em pullback} of $(H,\J)$ by $g$. Recall that $H^g$ depends
on the choice of a Hamiltonian which generates $g$. By Lemma
\ref{th:hg}, $\alpha_{H^g} = g^*\alpha_H$; therefore
\begin{equation} \label{eq:g:action}
a_{H^g} = \tilde{g}^*a_H + \mathrm{(constant)}
\end{equation}
and $\Crit(a_H) = \tilde{g}(\Crit(a_{H^g}))$ for every lift
$\tilde{g}: \tloops \longrightarrow \tloops$ of the action of $g$.

\begin{lemma} \label{th:functoriality}
For all $c_-,c_+ \in \Crit(a_{H^g})$, there
is a bijective map \[
\moduli(c_-,c_+;H^g,\J^g)/\R \longrightarrow
\moduli(\tilde{g}(c_-),\tilde{g}(c_+);H,\J)/\R.
\] Moreover, if $(H,\J)$ is a regular pair, 
$(H^g,\J^g)$ is also regular. \end{lemma}

\proof Let $\tilde{u}: \R \longrightarrow \tloops$ be a smooth path
whose projection to $\loops$ is given by $u \in \smooth(\R \times S^1,M)$.
Recall that $u$ satisfies \eqref{eq:floer} iff
\begin{equation} \label{eq:floer2}
\frac{d\tilde{u}}{ds} + \grad_{\J}a_H(\tilde{u}(s)) = 0
\end{equation}
Let $\tilde{v}(s) = \tilde{g}^{-1}(\tilde{u}(s))$. By \eqref{eq:g:metric}
and \eqref{eq:g:action}, $\tilde{u}$ satisfies \eqref{eq:floer2}
iff $\tilde{v}$ is a solution of
\begin{equation} \label{eq:floer3}
\frac{d\tilde{v}}{ds} + \grad_{\J^g}a_{H^g}(\tilde{v}(s)) = 0.
\end{equation}
The map of solution spaces defined in this way is clearly bijective
and $\R$-equivariant. Assuming that $(H,\J)$ is regular, we now
check that its pullback $(H^g,\J^g)$ satisfies the conditions for
a regular pair. (i) Let $c = [x,v]$ be a critical point of $a_{H^g}$,
$c' = \tilde{g}(c)$ and $\xi \in \smooth(x^*TM)$ with
$D^2a_{H^g}(c)(\xi,\cdot) = 0$. By \eqref{eq:g:action}, 
$D^2a_H(c')(Dg(x)\xi,\cdot) = 0$. Since
$c'$ is nondegenerate, $\xi = 0$, and therefore
$c$ is also nondegenerate. The last sentence in Definition
\ref{def:regular-pair}(i) follows from \eqref{eq:g:v}. (ii) can be 
proved using the same method as in Lemma \ref{th:naturality-of-regularity},
which also shows that a solution of \eqref{eq:floer3} and 
the corresponding solution of \eqref{eq:floer2} have the same index. 
Together with \eqref{eq:g:v}, this implies that $(H^g,\J^g)$ 
satisfies condition \ref{def:regular-pair}(iii). \qed

\begin{prop} \label{th:shift}
If $c \in \Crit(a_{H^g})$ is nondegenerate,
$\mu_H(\tilde{g}(c)) = \mu_{H^g}(c) - 2I(g,\tilde{g})$. \end{prop}

Recall that $\mu_{H^g}(c)$ is defined in the following way
(see \cite{salamon-zehnder92}, \cite{hofer-salamon95}): choose
a representative $(v,x)$ of $c$ and a symplectic trivialization
$\tau_c: x^*TM \longrightarrow S^1 \times (\R^{2n},\o_0)$ which 
can be extended over $v^*TM$. Let $\Psi_{H^g,c}: [0;1] \longrightarrow
\mySp(2n,\R)$ be the path given by
\[
\Psi_{H^g,c}(t) = \tau_c(t) D\phi_{H^g}^t(x(0)) \tau_c(0)^{-1},
\]
where $(\phi_{H^g}^t)_{t \in \R}$ is the Hamiltonian flow of $H^g$.
$\Psi_{H^g,c}(0) = \id$, and since $c$ is nondegenerate, 
$\mathrm{det}(\id - \Psi_{H^g,c}(1)) \neq 0$. A path $\Psi$ with these
properties has an index $\mu_1(\Psi) \in \Z$, and $\mu_{H^g}(c)$ is defined
by $\mu_{H^g}(c) = \mu_1(\Psi_{H^g,c})$. $\mu_1$ has the following property
\cite[Proposition 5]{floer-hofer93}:

\begin{lemma} \label{th:rescale} If $\Psi,\Psi'$ are related by
$\Psi'(t) = l(t) \Psi(t) l(0)^{-1}$ for some $l \in
\smooth(S^1,\mySp(2n,\R))$, $\mu_1(\Psi') = \mu_1(\Psi) - 2 \deg(l)$. \qed
\end{lemma}

\proof[Proof of Proposition \ref{th:shift}]
Let $(v',x')$ be a representative of $\tilde{g}(c)$ and $\tau_{\tilde{g}(c)}$ 
a trivialization of $(x')^*TM$ which can be extended over $(v')^*TM$.
$\mu_H(\tilde{g}(c))$ is defined using the path 
\[
\Psi_{H,\tilde{g}(c)}(t) = \tau_{\tilde{g}(c)}(t) D\phi^t_{H}(x'(0)) 
\tau_{\tilde{g}(c)}(0)^{-1}.
\]
$H^g$ is defined in such a way that
$\phi_H^t = g_t \, \phi_{H^g}^t \, g_0^{-1}$. Therefore
$
\Psi_{H,\tilde{g}(c)}(t) = l(t) \Psi_{H^g,c}(t) l(0)^{-1},
$ 
where $l(t) = \tau_{\tilde{g}(c)}(t) Dg_t(x(t)) \tau_c(t)^{-1}$. $l$ is
a loop in $\mySp(2n;\R)$, and $\deg(l) = I(g,\tilde{g})$ by definition. 
Proposition \ref{th:shift} now follows from Lemma \ref{th:rescale}. \qed

\begin{defn} \label{def:induced-maps}
Let $(H,\J)$ be a regular pair, $(g,\tilde{g}) \in \group$
and $(H^g,\J^g)$ the pullback of $(H,\J)$ by $g$. 
For $k \in \Z$, define an isomorphism
\[
CF_k(g,\tilde{g};H,\J,H^g): CF_k(H^g) \longrightarrow
CF_{k-2I(g,\tilde{g})}(H)
\]
by $CF_k(g,\tilde{g};H,\J,H^g)(\gen{c}) = \gen{\tilde{g}(c)}$,
extended in the obvious way to infinite sums of the generators
$\gen{c}$. Because of \eqref{eq:g:action}, this respects the 
finiteness condition for the formal sums in $CF_*(H)$. \end{defn}

In view of the definition of the
differential \eqref{eq:differential}, Lemma \ref{th:functoriality}
says that $CF_*(g,\tilde{g};H,\J,H^g)$ is an isomorphism of
chain complexes. We denote the induced isomorphisms on Floer homology by
\[
HF_k(g,\tilde{g};H,\J,H^g): HF_k(H^g,\J^g) \longrightarrow
HF_{k-2I(g,\tilde{g})}(H,\J).
\]
Consider two regular pairs $(H^-,\J^-)$, $(H^+,\J^+)$, two
Hamiltonians $K^-_g$, $K^+_g$ which generate $g$ and the
pullbacks $((H^-)^g,(\J^-)^g)$, $((H^+)^g,(\J^+)^g)$ using
$K^-_g$ and $K^+_g$, respectively. Let $(H,\J)$ be a homotopy
between $(H^-,\J^-)$ and $(H^+,\J^+)$. 
Choose $\psi \in \smooth(\R,\R)$ with $\psi|(-\infty;-1] = 0$ and
$\psi|[1;\infty) = 1$, and define $H^g \in \smooth(\R \times S^1
\times M,\R)$ by
\[
H^g(s,t,y) = H(s,t,g_t(y)) - (1-\psi(s)) K_g^-(t,g_t(y)) 
- \psi(s) K_g^+(t,g_t(y)).
\]
Together with the family $\J^g$ given by $J^g_{s,t} = Dg_t^{-1}J_{s,t}Dg_t$,
this is a homotopy from $((H^-)^g,(\J^-)^g)$ to $((H^+)^g,(\J^+)^g)$.

\begin{lemma} \label{th:homotopy-is-natural}
For all $c_- \in \Crit(a_{(H^-)^g}), c_+ \in \Crit(a_{(H^+)^g})$,
there is a bijective map
\[
\moduli^{\Phi}(c_-,c_+;H^g,\J^g) \longrightarrow
\moduli^{\Phi}(\tilde{g}(c_-),\tilde{g}(c_+);H,\J). 
\]
Moreover,
if $(H,\J)$ is a regular homotopy, $(H^g,\J^g)$ is also regular. \end{lemma}

We omit the proof, which is similar to that of 
Lemma \ref{th:functoriality}. 

\begin{cor} If $(H,\J)$ is a regular homotopy, the diagram
\begin{equation*}
\xymatrix{
CF_*((H^-)^g,(\J^-)^g) \ar[rrrr]^{CF_*(g,\tilde{g};H^-,\J^-,(H^-)^g)} 
\ar[d]^{\Phi_*(H^g,\J^g)}
&&&& CF_*(H^-,\J^-) \ar[d]^{\Phi_*(H,\J)}\\
CF_*((H^+)^g,(\J^+)^g) \ar[rrrr]^{CF_*(g,\tilde{g};H^+,\J^+,(H^+)^g)} 
&&&& CF_*(H^+,\J^+)
} 
\end{equation*} 
commutes. \end{cor}

This is proved by a straightforward computation using Lemma
\ref{th:homotopy-is-natural}. It follows that for every 
$\pair \in \group$, there is a unique automorphism
\[
HF_*\pair: \floer \longrightarrow \floer
\]
independent of the choice of a regular pair.
We list some properties of these maps which are clear
from the definition.

\begin{proposition} \label{prop:further-properties}
\begin{enumerate}
\item
$HF_*\pair$ is an automorphism of $\floer$ as a
$\novi$-module.
\item For $\pair = \id_{\group}$, $HF_*\pair = \id_{\floer}$.
\item \label{item:lambda-multiplication}
If $\pair = (\id,\gamma)$ for some $\gamma \in \Gamma$, $HF_*\pair$ 
is equal to the multiplication by $\gen{\gamma} \in \novi$.
\item \label{item:functor}
For $(g_1,\tilde{g}_1), (g_2,\tilde{g}_2) \in \group$,
$
HF_*(g_1 g_2,\tilde{g}_1 \tilde{g}_2) = 
HF_*(g_1,\tilde{g}_1) HF_*(g_2,\tilde{g}_2)
$. \qed
\end{enumerate}
\end{proposition}

%
\newcommand{\modh}[2]{\moduli^h(#1,#2;\bar{H},\bar{\J})}
\newcommand{\cmodh}[2]{\overline{\moduli}^h(#1,#2;\bar{H},\bar{\J})}

\section{Homotopy invariance \label{sec:homotopy}}

Let $(g_{r,t})_{0 \leq r \leq 1, t \in S^1}$ be a smooth family of Hamiltonian
automorphisms of $\mo$ with $g_{r,0} = \id$ for all $r$. This family
defines a path $(g_r)_{0 \leq r \leq 1}$ in $G$. Let 
$(g_r,\tilde{g}_r)_{0 \leq r \leq 1}$ be a smooth lift of this path to
$\group$. The aim of this section is to prove the following

\begin{prop} \label{th:homotopy}
For $(g_r,\tilde{g}_r)_{0 \leq r \leq 1}$ as above,
\[
HF_*(g_0,\tilde{g}_0) = HF_*(g_1,\tilde{g}_1):
HF_*\mo \longrightarrow HF_*\mo. 
\]\end{prop}

Choose a smooth family $(K_r)_{0 \leq r \leq 1}$ of Hamiltonians such 
that $K_r$ generates $g_r$. Let $(H,\J)$ be a regular pair and
$(H^{g_r},\J^{g_r})$ its pullbacks (using $K_r$). The maps induced
by $(g_0,\tilde{g}_0)$ and $(g_1,\tilde{g}_1)$ are
\begin{align*}
HF_*(g_0,\tilde{g}_0;H,\J,H^{g_0}): HF_*(H^{g_0},\J^{g_0}) &\longrightarrow 
HF_*(H,\J),\\
HF_*(g_1,\tilde{g}_1;H,\J,H^{g_1}): HF_*(H^{g_1},\J^{g_1}) &\longrightarrow 
HF_*(H,\J).
\end{align*}
Proposition \ref{th:homotopy} says that if $(H',\J')$ is a regular
homotopy from $(H^{g_1},\J^{g_1})$ to $(H^{g_0},\J^{g_0})$,
\begin{equation} \label{eq:htp-one}
HF_*(g_1,\tilde{g}_1;H,\J,H^{g_1}) = HF_*(g_0,\tilde{g}_0;H,\J,H^{g_0})
\Phi_*(H',\J').
\end{equation}
Because of Proposition \ref{prop:further-properties}\ref{item:functor}, 
it is sufficient to consider the case $(g_0,\tilde{g}_0) = {\id}_{\group}$.
Choose $K_0 = 0$. Then $(H^{g_0},\J^{g_0}) = (H,\J)$,
$HF_*(g_0,\tilde{g}_0; H,\J,H) = \id_{HF_*(H,\J)}$, and
\eqref{eq:htp-one} is reduced to
\begin{equation} \label{eq:htp-two}
\Phi_*(H',\J') HF_*(g_1,\tilde{g}_1;H,\J,H^{g_1})^{-1} = \id_{HF_*(H,\J)}.
\end{equation}
The rest of the section contains the proof of this equation.

\begin{defn} Let $(H,\J)$, $(H^{g_r},\J^{g_r})$ and
$(H',\J')$ be as above. A {\em deformation of homotopies} compatible with
them consists of a function 
$\bar{H} \in \smooth([0;1] \times \R \times S^1 \times
M,\R)$ and a family of almost complex structures $\bar{\J} =
(\bar{J}_{r,s,t}) \in \modJreg{[0;1] \times \R \times S^1}$, such that
\begin{align*}
& \bar{H}(r,s,t,y) = H^{g_r}(t,y), \; \bar{J}_{r,s,t} = J^{g_r}_t 
&& \text{ for } s \leq -1, \\
& \bar{H}(r,s,t,y) = H(t,y), \; \bar{J}_{r,s,t} = J_t
&& \text{ for } s \geq 1, \\
& \bar{H}(0,s,t,y) = H(t,y), \; \bar{J}_{0,s,t} = J_t && \text{ and }\\
& \bar{H}(1,s,t,y) = H'(s,t,y), \; \bar{J}_{1,s,t} = J_{s,t}'. &&
\end{align*} \end{defn}

Let $(\bar{H},\bar{\J})$ be a deformation of homotopies. Consider a
pair $(r,u) \in [0;1] \times \smooth(\R \times S^1,M)$ such that
\begin{equation} \label{eq:parametrized-floer}
\partials{u} + \bar{J}_{r,s,t}(u(s,t))\left(
\partialt{u} - X_{\bar{H}}(r,s,t,u(s,t))\right) = 0.
\end{equation}
We say that $u$ converges to
$c_-,c_+ \in \Crit(a_H)$ if there is a smooth path
$\tilde{u}: \R \longrightarrow \tloops$ with
\[
\lim_{s \rightarrow -\infty} \tilde{u}(s) = \tilde{g}_r^{-1}(c_-),
\lim_{s \rightarrow +\infty} \tilde{u}(s) = c_+
\]
such that $u(s,\cdot) = p(\tilde{u}(s))$ for all $s$.
We denote the space of such pairs $(r,u)$ with limits
$c_\pm$ by $\moduli^h(c_-,c_+;\bar{H},\bar{\J})$.
The linearization of \eqref{eq:parametrized-floer}
at a pair $(r,u) \in \moduli^h(c_-,c_+;\bar{H},\bar{\J})$ 
is given by an operator
\[
D_{\bar{H},\bar{\J}}^h(r,u): \R \times \smooth(u^*TM)
\longrightarrow \smooth(u^*TM).
\]
Since $c_-,c_+$ are nondegenerate critical points of $a_H$,
$\tilde{g}_r^{-1}(c_-)$ is a nondegenerate critical point
of $a_{H^{g_r}}$ for all $r$. It follows that the Sobolev
completion $(p>2)$
\begin{equation} \label{eq:h-operator}
D_{\bar{H},\bar{\J}}^h(r,u): \R \times W^{1,p}(u^*TM)
\longrightarrow L^p(u^*TM)
\end{equation}
is a Fredholm operator of index
$\mu_{H^{g_r}}(\tilde{g}_r^{-1}(c_-)) - \mu_H(c_+) + 1$. 
By Lemma \ref{th:shift}, $\mu_H(c_-) =
\mu_{H^{g_r}}(\tilde{g}_r^{-1}(c_-)) - 2 I(g_r,\tilde{g}_r)$.
However, since $(g_r,\tilde{g}_r)$ is homotopic to the identity
in $\group$, $I(g_r,\tilde{g}_r) = 0$. It follows that
$\ind D_{\bar{H},\bar{\J}}^h(r,u) = \mu_H(c_-) - \mu_H(c_+) + 1$.

$(\bar{H},\bar{\J})$ is called regular if \eqref{eq:h-operator} 
is onto for all $(r,u)$, and
$
(r,s,t,u(s,t)) \notin V_0(\bar{\J})
$
for all $(r,s,t) \in [0;1] \times \R \times S^1$ and
$u \in \moduli^h(c_-,c_+;\bar{H},\bar{\J})$ such that
$\mu_H(c_-) \leq \mu_H(c_+)$. An analogue of Theorem
\ref{th:transv-floer} ensures that regular deformations
of homotopies exists. Regularity implies that the spaces
$\moduli^h(c_-,c_+;\bar{H},\bar{\J})$ are smooth manifolds.
The boundary $\partial\moduli^h(c_-,c_+;\bar{H},\bar{\J})$ consists
of solutions of \eqref{eq:parametrized-floer} with
$r = 0$ or $1$. For $r = 0$, \eqref{eq:parametrized-floer} is
\begin{equation} \label{eq:r-is-zero}
\partials{u} + J_t(u)\left(\partialt{u} - X_H(t,u)\right) = 0
\end{equation}
and for $r = 1$, it is 
\begin{equation} \label{eq:r-is-one}
\partials{u} + J_{s,t}'(u)\left(\partialt{u} - X_{H'}(s,t,u)\right) = 0.
\end{equation}
Since $(H,\J)$ is a regular pair and $(H',\J')$ a regular homotopy,
this implies that $\partial\moduli^h(c_-,c_+;H,\J) = \emptyset$ if
$\mu_H(c_+) = \mu_H(c_-) + 1$.

\begin{lemma} \label{th:parametrized-compactness}
Let $(\bar{H},\bar{\J})$ be a regular deformation of homotopies
and $c_- \in \Crit_{k_-}(a_H)$, $c_+ \in \Crit_{k_+}(a_H)$. 
\begin{enumerate}
\item \label{item:zero-d}
If $k_+ = k_- + 1$, $\modh{c_-}{c_+}$ is a finite set. 
\item If $k_+ = k_-$, $\modh{c_-}{c_+}$ is one-dimensional,
and there is a smooth compactification $\cmodh{c_-}{c_+}$ whose boundary
$\partial\cmodh{c_-}{c_+}$ is the disjoint union of 
$\partial\modh{c_-}{c_+}$,
\begin{equation} \label{eq:comp-component-1}
\modh{c_-}{c} \times (\moduli(c,c_+;H,\J)/\R)
\end{equation}
for all $c \in \Crit_{k_+ + 1}(a_H)$ and
\begin{equation} \label{eq:comp-component-2}
(\moduli(c_-,c';H,\J)/\R) \times \modh{c'}{c_+} 
\end{equation}
for all $c' \in \Crit_{k_+ - 1}(a_H)$. \end{enumerate}
\end{lemma}

This compactification is constructed in the same way as that used in
\cite[Lemma 6.3]{salamon-zehnder92}. We omit the proof. 
The fact that the spaces $\modh{c_-}{c_+}$ of 
dimension $\leq 1$ can be compactified without including limit points
with `bubbles' is a consequence of the regularity
condition (compare \cite[Theorem 5.2]{hofer-salamon95}).
The limits \eqref{eq:comp-component-2} arise in the following way:
let $(r_m,u_m)$ be a sequence in $\modh{c_-}{c_+}$ with
$\lim_m r_m = r$ and such that $u_m$ converges uniformly on compact
subsets to some $u_\infty \in \smooth(\R \times S^1,M)$. Assume that
there are $s_m \in \R$, $\lim_m s_m = \infty$, such that the
translates $\hat{u}_m(s,t) = u_m(s-s_m,t)$
converge on compact subsets to a map $\hat{u}_\infty$ with
$\partial\hat{u}_\infty/\partial s \not\equiv 0$.
$\hat{u}_\infty$ is a solution of
\[
\partials{\hat{u}_\infty} + J_t^{g_r}(\hat{u}_\infty)\left(
\partialt{\hat{u}_\infty} - X_{H^{g_r}}(t,\hat{u}_\infty)\right) = 0.
\]
Therefore $\hat{v}_\infty(s,t) = g_{r,t}(\hat{u}_\infty(s,t))$ defines
a solution of \eqref{eq:r-is-zero}. $(\hat{v}_\infty,(r,u_\infty))$ is
the limit of the sequence $u_m$ in \eqref{eq:comp-component-2}.
The points \eqref{eq:comp-component-1} in the compactification
arise by a similar, but simpler process using translations with
$s_m \longrightarrow -\infty$.

\begin{defn} \label{def:chain-homotopy} Let $(\bar{H},\bar{\J})$ be
a regular deformation of homotopies. For $k \in \Z$, define
\begin{align*}
h_k(\bar{H},\bar{\J}): CF_k(H) &\longrightarrow CF_{k+1}(H),\\
h_k(\bar{H},\bar{\J})(\gen{c_-}) &=
\sum_{c_+} \#\modh{c_-}{c_+} \gen{c_+}.
\end{align*} \end{defn}
To show that this formal sum is indeed an element of $CF_*(H,\J)$,
a slightly stronger version of Lemma
\ref{th:parametrized-compactness}\ref{item:zero-d}
is necessary; we omit the details.

\begin{lemma} \label{th:it-is-a-chain-homotopy} For all $k$,
\begin{multline*}
\partial_{k+1}(H,\J) h_k(\bar{H},\bar{\J}) + h_{k-1}(\bar{H},\bar{\J})
\partial_k(H,\J) =\\
= \Phi_k(H',\J')CF_k(g_1,\tilde{g}_1;H,\J,H^{g_1})^{-1} - \id_{CF_k(H)}. 
\end{multline*}
\end{lemma} 

\proof By definition,
\begin{multline*}
\Phi_k(H',\J') CF_k(g_1,\tilde{g}_1;H,\J,H^{g_1})^{-1} \gen{c_-} =\\
= \sum_{c_+} \#\moduli^{\Phi}(\tilde{g}_1^{-1}(c_-),c_+;H',\J')\gen{c_+}.
\end{multline*}
Comparing \eqref{eq:r-is-one} with \eqref{eq:continuation},
one sees that $u \in \moduli^\Phi(\tilde{g}_1^{-1}(c_-),c_+;H',\J')$
iff $(1,u) \in \modh{c_-}{c_+}$. The other part of
$\partial\modh{c_-}{c_+}$ is given by solutions of \eqref{eq:r-is-zero}
with limits $c_\pm$. If $c_-,c_+$ have the same Conley-Zehnder index,
any regular solution of \eqref{eq:r-is-zero} is stationary in the sense
that $\partial u/\partial s = 0$, and therefore
this part of $\partial\modh{c_-}{c_+}$ is empty unless
$c_- = c_+$, in which case it contains a single point. It follows that
\begin{multline*}
\left(\Phi_k(H',\J')CF_k(g_1,\tilde{g}_1;H,\J,H^{g_1})^{-1} - \id\right)
\gen{c_-} = \\
= \sum_{c_+} \#\partial\modh{c_-}{c_+} \gen{c_+}
\end{multline*}
(the sign is irrelevant since we use $\Z/2$-coefficients).
Similarly, the coefficients of $\gen{c_+}$ in 
$\partial_{k+1}(H,\J) h_k(\bar{H},\bar{\J})\gen{c_-}$ and 
$h_{k-1}(\bar{H},\bar{\J}) \partial_k(H,\J)\gen{c_-}$ are 
given by the number of points in
\eqref{eq:comp-component-1} and \eqref{eq:comp-component-2}. 
The statement now follows from the
fact that $\partial\cmodh{c_-}{c_+}$ contains an even number of points. \qed

Lemma \ref{th:it-is-a-chain-homotopy} shows that
$\Phi_*(H',\J') CF_*(g_1,\tilde{g}_1;H,\J,H^{g_1})^{-1}$
is chain homotopic to the identity, which implies
\eqref{eq:htp-two}.

\newcommand{\modpp}{\moduli^{\mathrm{PP}}(c_-,c_0,c_+;H,\J)}
\newcommand{\gmodpp}{\moduli^{\mathrm{PP}}(c_-,c_0,c_+;H^g,\J^g)}
\newcommand{\ggmodpp}{\moduli^{\mathrm{PP}}(\tilde{g}(c_-),c_0,
\tilde{g}(c_+);H,\J)}
\newcommand{\perturb}[6]{\mathcal{P}(#1,#2,#3,#4,#5,#6)}
\newcommand{\perturbation}{\perturb{H^-}{\J^-}{H^0}{\J^0}{H^+}{\J^+}}
\newcommand{\onefourth}{\frac{1}{4}}
\section{Pair-of-pants product \label{sec:pair-of-pants}}

The main result of this section, Proposition \ref{th:pair-of-pants-linear},
describes the behaviour of the pair-of-pants product on Floer homology
under the maps $HF_*(g,\tilde{g})$. It turns out that $HF_*\pair$ does
not preserve the ring structure; rather, it is an automorphism of
$\floer$ as a module over itself.

The definition of the pair-of-pants product for weakly monotone symplectic
manifolds can be found in {\pss} or \cite[Chapter 10]{mcduff-salamon}
(a detailed construction in the case $[\o]|\pi_2(M) = 0$ is given in
\cite{schwarz95}). We need to modify this definition slightly by enlarging
the class of admissible perturbations. This would be necessary in any case
to adapt it to the use of time-dependent almost complex structures; 
a further enlargement ensures that the relevant moduli spaces
transform nicely under (a subset of) $\group$. 
We now sketch this modified definition. 

Consider the punctured surface
$\Sigma = \R \times S^1 \setminus (0;0)$. It has a third tubular end
\[
e: \R^- \times S^1 \longrightarrow \Sigma, \quad
e(s,t) = (\frac{1}{4}e^{2\pi s}\cos(2\pi t), 
\frac{1}{4}e^{2\pi s}\sin(2\pi t))
\]
around the puncture, and the surface $\widehat{\Sigma}$ obtained by
capping off this end with a disc can be identified with $\R \times S^1$.
For $u \in \smooth(\Sigma,M)$ we will denote $u \circ e$ by $u^e$.
Two maps $u \in \smooth(\Sigma,M)$ and $v_0 \in \smooth(D^2,M)$ which
satisfy $\lim_{s \rightarrow -\infty} u^e(s,\cdot) = v_0|\partial D^2$
can be glued together to a map $u \glue v_0: \widehat{\Sigma}
\longrightarrow M$. By identifying $\widehat{\Sigma} = \R \times S^1$,
this gives a path $\R \longrightarrow \Lambda M$. On the two 
remaining ends, $u \glue v_0$ has the same asymptotic behaviour as $u$.
In particular, if there are $x_-, x_+ \in \loops$ such that
$\lim_{s \rightarrow \pm\infty} u(s,\cdot) = x_\pm$, the path
$u \glue v_0$ lies in $\loops$ and converges to $x_-$, $x_+$.

\begin{defn} We say that $u \in \smooth(\Sigma,M)$ {\em converges
to} $c_-,c_0,c_+ \in \tloops$ if
\[
\lim_{s \rightarrow \pm\infty} u(s,\cdot) = p(c_\pm),
\quad \lim_{s \rightarrow -\infty} u^e(s,\cdot) = p(c_0)
\]
and if for $c_0 = [v_0,x_0]$, the path $u \glue v_0: \R \longrightarrow
\loops$ has a lift $\widetilde{u \glue v_0}: \R \longrightarrow \tloops$
with limits $c_-,c_+$. \end{defn}

This is independent of the choice of the representative $(v_0,x_0)$ for
$c_0$ and of the choices involved in the gluing.

Let $(H^-,\J^-)$, $(H^0,\J^0)$ and $(H^+,\J^+)$ be regular pairs.
Choose $\J = \! (J_z)_{z \in \Sigma} \in \modJreg{\Sigma}$ and
$H \in \smooth(\Sigma \times M,\R)$ with $H|e([-1;0] \times S^1) \times M
= 0$, such that $(H,\J)$ agrees with one regular pair over each end
of $\Sigma$, more precisely:
\begin{align*}
& H(-s,t,y) = H^-(t,y), \; J_{(-s,t)} = J_t^-, \\
& H(s,t,y) = H^+(t,y), \; J_{(s,t)} = J_t^+, \\
& H(e(-s,t),y) = H^0(t,y), \; J_{e(-s,t)} = J_t^0 
\end{align*}
for $s \geq 2$, $t \in S^1$ and $y \in M$. The space of such
$(H,\J)$ will be denoted by $\perturb{H^-}{\J^-}{H^0}{\J^0}{H^+}{\J^+}$.

$u \in \smooth(\Sigma,M)$ is called $(H,\J)$-holomorphic if
\[
\partials{u}(s,t) + J_{s,t}(u(s,t))\left( 
\partialt{u}(s,t) - X_H(s,t,u(s,t))\right)
= 0
\]
for $(s,t) \in \Sigma\setminus e((-\infty;-1] \times S^1)$ and
\[
\partials{u^e}(s,t) + J_{e(s,t)}(u^e(s,t)) \left(
\partialt{u^e}(s,t) - X_H(e(s,t),u^e(s,t))\right) = 0
\]
for $(s,t) \in \R^- \times S^1$. 
For $z = (s,t) \in e([-1;0] \times S^1)$, the first equation be written
as $\bar{\partial}_{J_z}u(z) = 0$ using the complex structure
on $\Sigma \iso (\C/i\Z) \setminus \! (0;0)$. Similarly, the second equation
is $\bar{\partial}_{J_{e(s,t)}}u^e(s,t) = 0$ 
for $(s,t) \in [-1;0] \times S^1$. Since $e$ is holomorphic, this
implies that the two equations match up smoothly.

\begin{defn} For $c_- \in \Crit(a_{H^-}), c_0 \in \Crit(a_{H^0})$
and $c_+ \in \Crit(a_{H^+})$, $\modpp$ is the space of
$(H,\J)$-holomorphic maps which converge to $c_-,c_0,c_+$
in the sense discussed above. \end{defn}

The linearization of the two equations for an $(H,\J)$-holomorphic map at
$u \in \modpp$ is given by single differential operator
$D^\Sigma_{H,\J}(u)$ on $\Sigma$ which becomes a Fredholm operator
in suitable Sobolev spaces. Following the same method as in 
the previous sections, we call $(H,\J)$ regular if all these 
operators are onto and if for all $(H,\J)$-holomorphic maps $u$ with 
$\ind(D^\Sigma_{H,\J}(u)) \leq 1$ and all $z \in \Sigma$, 
$(z,u(z)) \notin V_0(\J)$. A transversality result analogous to Theorem 
\ref{th:transv-floer} shows that the set of regular $(H,\J)$ is dense in 
$\perturb{H^-}{\J^-}{H^0}{\J^0}{H^+}{\J^+}$. If $(H,\J)$ is
regular and $\mu_{H^+}(c_+) = \mu_{H^-}(c_-) + \mu_{H^0}(c_0) - 2n$, 
the space $\modpp$ is finite. The pair-of-pants product
\[
PP_{i,j}(H,\J): CF_{i}(H^-) \otimes CF_{j}(H^0)
\longrightarrow CF_{i+j-2n}(H^+)
\]
for regular $(H,\J)$ is defined by
\[
PP_{i,j}(H,\J)(\gen{c_-} \otimes \gen{c_0}) = \sum_{c_+} \#\modpp \gen{c_+}.
\]
A compactness theorem shows that this formal sum lies in $CF_*(H^+)$.
The grading is a consequence of our convention for
the Conley-Zehnder index.

The construction of the product is completed in the following
steps: using the spaces $\moduli^{\mathrm{PP}}$ of dimension $1$ and their
compactifications, one shows that $PP_*(H,\J)$ is a chain homomorphism. 
The induced maps
\begin{equation} \label{eq:ppp}
PP_*(H,\J): HF_*(H^-,\J^-) \times HF_*(H^0,\J^0) \longrightarrow
HF_*(H^+,\J^+)
\end{equation}
are independent of the choice of $(H,\J)$ (the starting point for
the proof is that $\perturb{H^-}{\J^-}{H^0}{\J^0}{H^+}{\J^+}$ 
is contractible). Finally, a gluing argument proves that the 
maps \eqref{eq:ppp} for different regular pairs are related by 
continuation isomorphisms. Therefore they define a unique product
\[
\pairofpants: HF_*\mo \times HF_*\mo \longrightarrow HF_*\mo.
\]
Using the assumption {\wplus}, the proofs of these properties 
for the product as defined in {\pss} can be easily adapted to 
our slightly different setup.

\begin{prop} \label{th:pair-of-pants-linear}
For all $(g,\tilde{g}) \in \group$ and $a,b \in HF_*\mo$,
\[
HF_*(g,\tilde{g})(a \pairofpants b) =
HF_*(g,\tilde{g})(a) \pairofpants b.
\]
\end{prop}

Because of the homotopy invariance of $HF_*(g,\tilde{g})$
(Proposition \ref{th:homotopy}) it is sufficient to consider
the case where $g_t = \id_M$ for $t \in [-\onefourth;\onefourth] \subset S^1$
(clearly, any path component of $\group$ contains a $(g,\tilde{g})$
with this property). Choose a Hamiltonian $K_g$ which generates $g$
such that $K_g(t,\cdot) = 0$ for $t \in [-\onefourth;\onefourth]$.
Let $((H^{\pm})^g,(\J^{\pm})^g)$ be the pullback of $(H^\pm,\J^\pm)$ using
$K_g$. The proof of Proposition \ref{th:pair-of-pants-linear}
relies on the following analogue of the `pullback' of a regular
pair: for $(H,\J) \in \perturbation$, define $(H^g,\J^g)$ by
\[
J^g_{(s,t)} = Dg_t^{-1} J_{(s,t)} Dg_t, \quad
H^g(s,t,y) = H(s,t,g_t(y)) - K_g(t,g_t(y))
\]
for $(s,t) \in \Sigma$, $y \in M$. Because $g_t = \id_M$ for 
$t \in [-\onefourth;\onefourth]$ and
\[
\mathrm{im}(e) \subset (\R \times [-1/4;1/4]) \setminus (0;0) 
\subset \Sigma, 
\] 
$(H^g,\J^g)$ satisfies 
$H^g(e(s,t),\cdot) = H(e(s,t),\cdot)$ and $J^g_{e(s,t)} = J_{e(s,t)}$ for
all $(s,t) \in \R^- \times S^1$. In particular, $H^g$ vanishes on
$e([-1;0] \times S^1) \times M$. It follows that
$(H^g,\J^g) \in \perturb{(H^-)^g}{(\J^-)^g}{H^0}{\J^0}{(H^+)^g}{(\J^+)^g}$.
The next Lemma is the analogue of Lemma \ref{th:functoriality}:

\begin{lemma} \label{th:pp-functoriality} For all 
$(c_-,c_0,c_+) \in \Crit(a_{(H^-)^g}) \times \Crit(a_{H^0})
\times \Crit(a_{(H^+)^g})$, there is a bijective map
\[
\gmodpp \longrightarrow \ggmodpp.
\]
Moreover, if $(H,\J)$ is regular, so is $(H^g,\J^g)$.
\end{lemma}

\proof A straightforward computation shows that
if $u,v \in \smooth(\Sigma,M)$ are related by $u(s,t) = g_t(v(s,t))$, 
$u$ is $(H,\J)$-holomorphic iff $v$ is $(H^g,\J^g)$-holomorphic 
(note that $v^e = u^e$, and that the second equation
is the same in both cases). Now assume that $v$ converges
to $c_-,c_0,c_+$ as defined above, and choose a 
representative $(v_0,x_0)$ of $c_0$. The maps
$u \glue v_0, v \glue v_0: \R \times S^1 \longrightarrow M$ 
can be chosen such that
\[
(u \glue v_0)(s,t) = g_t((v \glue v_0)(s,t)).
\]
By assumption, the path $\R \longrightarrow \loops$ given by
$v \glue v_0$ has a lift $\widetilde{v \glue v_0}: \R \longrightarrow
\tloops$ with limits $c_-,c_+$. Define
\[
\widetilde{u \glue v_0}(s) = \tilde{g}(\widetilde{v \glue v_0}(s)).
\]
Clearly, $\widetilde{u \glue v_0}$ is a lift of $u \glue v_0$ with
limits $\tilde{g}(c_-)$, $\tilde{g}(c_+)$. Therefore $u$ converges
to $\tilde{g}(c_-), c_0, \tilde{g}(c_+)$. The proof of regularity
is similar to that of Lemma \ref{th:naturality-of-regularity}; we omit
the details. \qed

From Lemma \ref{th:pp-functoriality} and the definition
of $PP_*(H,\J)$, it follows that
\begin{multline*}
CF_*(g,\tilde{g};H^+,\J^+,(H^+)^g) PP_*(H^g,\J^g)(\gen{c_-}
\otimes \gen{c_0}) =\\
 = PP_*(H,\J)(\gen{\tilde{g}(c_-)} \otimes \gen{c_0})
\end{multline*}
for every regular $(H,\J) \in \perturbation$,
$c_- \in \Crit(a_{(H^-)^g})$ and $c_0 \in \Crit(a_{H^0})$.
Proposition \ref{th:pair-of-pants-linear} follows directly
from this.

\newcommand{\modeo}{\mathcal{J}(E,\O)}
\newcommand{\modeoreg}{\mathcal{J}^{\mathrm{reg}}(E,\O)}
\newcommand{\modeoregz}{\mathcal{J}^{\mathrm{reg},z_0}(E,\O)}
\newcommand{\cu}{\mathcal{C}}
\newcommand{\cusp}{\cu_{r,k}(J)}
\newcommand{\modhat}{\widehat{\mathcal{J}}(j,\J)}
\newcommand{\modhatreg}{\widehat{\mathcal{J}}^{\mathrm{reg}}(j,\J)}
\newcommand{\jjsec}{\mathcal{S}(j,\hatj)}
\newcommand{\duniv}{D^{\mathrm{univ}}}
\newcommand{\U}{\mathcal{U}}

\section{Pseudoholomorphic sections \label{sec:holomorphic-sections}}

Let $(E,\O)$ be a symplectic fibre bundle over $S^2$ with fibre
$\mo$ and $\pi: E \longrightarrow S^2$ the projection. 
We will denote by $\modeo$ the space of families 
$\J = (J_z)_{z \in S^2}$ of almost complex
structures on the fibres of $E$ such that $J_z$ is $\O_z$-compatible
for all $z$. For $\J \in \modeo$ and $k \in \Z$, let $\moduli_k^s(\J)$ 
be the space of pairs $(z,w) \in S^2 \times \smooth(\CP{1},E)$ such that 
$w$ is a simple $J_z$-holomorphic curve in $E_z$ with 
$c_1(TE_z,\O_z)(w) = k$. Because $(E,\O)$ is locally trivial, these spaces 
have the same properties as the spaces of holomorphic curves in $M$ with 
respect to a two-parameter family of almost complex structures (see
section \ref{sec:floer-homology}). In particular, there is a dense subset
$\modeoreg \subset \modeo$ such that\footnote{This uses the
assumption {\wplus}; again, weak monotonicity would not be sufficient.}
for $\J \in \modeoreg$, $\moduli_k^s(\J) = \emptyset$ for all $k<0$
and $\moduli_0^s(\J)$ is a manifold of dimension
$\dim \moduli_0^s(\J) = \dim E$. It follows that the image of the
evaluation map
\[
\eta: \moduli_0^s(\J) \times_{PSL(2,\C)} \CP{1} \longrightarrow E
\]
is a subset of codimension $4$.

Following the convention of section \ref{sec:elementary}, $(E,\O)$
is equipped with a preferred isomorphism $i: \mo \longrightarrow
(E_{z_0},\O_{z_0})$ for the marked point $z_0 \in S^2$. We need
to recall the transversality theory of cusp-curves; a reference
is \cite[Chapter 6]{mcduff-salamon}.

\begin{defn} Let $J$ be an $\o$-compatible almost complex structure
on $M$. A simple $J$-holomorphic cusp-curve with $r \geq 1$ components is a
\[
v = (w_1,\dots w_r, t_1,\dots t_r, t_1',\dots t_r')
\in \smooth(\CP{1},M)^r \times (\CP{1})^{2r}
\]
with the following properties:
\begin{enumerate} \item For $i = 1 \dots r$,
$w_i$ is a simple $J$-holomorphic curve;
\item $\im(w_i) \neq \im(w_j)$ for $i \neq j$;
\item $w_i(t_i) = w_{i+1}(t_{i+1}')$ for $i = 1 \dots r-1$.
\end{enumerate} \end{defn}

The Chern number of $v$ is defined by $c_1(v) = \sum_i c_1(w_i)$.
$PSL(2,\C)^r$ acts freely on the space of simple cusp-curves 
with $r$ components by 
\begin{multline*}
(a_1 \dots a_r) \cdot (w_1,\dots w_r, t_1,\dots t_r, t_1',\dots t_r') =\\
(w_1 \circ a_1^{-1},\dots w_r \circ a_r^{-1}, a_1(t_1),\dots
a_r(t_r), a_1(t_1'),\dots a_r(t_r')).
\end{multline*}
Let $\cusp$ be the quotient of the space of cusp-curves with $r$
components and Chern number $k$ by this action, and
\[
\eta_1, \eta_2: \cusp \longrightarrow M
\]
the maps given by $\eta_1(v) = w_1(t_1')$, $\eta_2(v) = w_r(t_r)$ for
$v$ as above. For a generic $J$, $\cusp$ is a smooth manifold of
dimension $2n + 2k - 2r$ \cite[Theorem 5.2.1(ii)]{mcduff-salamon}.
Let $\modeoregz \subset \modeoreg$ be the subset of families
$\J = (J_z)_{z \in S^2}$ such that $J = Di^{-1} J_{z_0} Di$ has this
regularity property; $\modeoregz$ is dense in $\modeo$.

Let $j$ be a positively oriented complex structure on $S^2$ and
$\J \in \modeo$. We call an almost complex structure $\jhat$ on $E$
compatible with $j$ and $\J$ if $D\pi \circ \hatj = j \circ D\pi$ and 
$\hatj|E_z = J_z$ for all $z \in S^2$. The space of such $\hatj$ will be 
denoted by $\modhat$. If we fix a $\hatj_0 \in \modhat$, any other
such $\hatj$ is of the form $\hatj_0 + \theta \circ D\pi$,
where $\theta: \pi^*(TS^2,j) \longrightarrow (TE^v,\J)$
is a smooth $\C$-antilinear vector bundle homomorphism
(recall that $TE^v = \ker D\pi$ is the tangent bundle along the fibres); 
conversely, $\hatj_0 + \theta \circ D\pi \in \modhat$ for any such $\theta$.

For $j,\J$ as above and $\hatj \in \modhat$, we call a smooth section
$s: S^2 \longrightarrow E$ of $\pi$ {\em$(j,\hatj)$-holomorphic} if
\begin{equation} \label{eq:holomorphic-section}
ds \circ j = \jhat \circ ds.
\end{equation}

If $E = S^2 \times M$ and $J_z$ is independent of
$z$, this is the `inhomogeneous Cauchy-Riemann equation' of 
\cite{ruan-tian94} for maps $S^2 \longrightarrow M$ 
with an inhomogeneous term determined by the choice of $\hatj$. 

Let $\jjsec$ be the space of $(j,\jhat)$-holomorphic sections of $E$.
For every section $s$, $\J$ induces an almost complex structure
$s^*\J$ on $s^*TE^v$, given by $(s^*\J)_z = (J_z)_{s(z)}$. 
The linearization of \eqref{eq:holomorphic-section} at $s \in \jjsec$
is a differential operator
\[
D_{\jhat}(s): \smooth(s^*TE^v) \longrightarrow \Omega^{0,1}(s^*TE^v,s^*\J)
\]
on $S^2$. $D_{\jhat}(s)$ differs  
from the $\bar{\partial}$-operator of $(s^*TE^v,s^*\J)$
by a term of order zero. Therefore it is a Fredholm operator with index
$d(s) = 2n + 2c_1(TE^v,\O)(s)$. If $D_{\jhat}(s)$ is onto,
$\jjsec$ is a smooth $d(s)$-dimensional manifold near $s$.

\begin{defn} Let $j$ be a positively oriented complex structure on $S^2$,
$\J = (J_z)_{z \in S^2} \in \modeoregz$ and $J$ the almost complex 
structure on $M$ given by $J = Di^{-1} J_{z_0} Di$. 
$\hatj \in \modhat$ is called regular if
$D_{\jhat}(s)$ is onto for all $s \in \jjsec$, the map
\[
\ev: S^2 \times \jjsec \longrightarrow E, \; \ev(z,s) = s(z)
\]
is transverse to $\eta$ and the map
\[
\ev_{z_0}: \jjsec \longrightarrow M, \; \ev_{z_0}(s) = i^{-1}(s(z_0))
\]
is transverse to $\eta_1$. \end{defn}

\begin{prop} The subspace $\modhatreg \subset \modhat$ 
of regular $\jhat$ is $\smooth$-dense.
\end{prop}

\proof[Proof (sketch)]
The strategy of the proof is familiar. The tangent space of 
$\modhat$ at any point is a subspace of the space of bundle homomorphisms
$I: TE \longrightarrow TE$. As explained above, it contains precisely those
$I$ which can be written as $I = \theta \circ D\pi$, where 
$\theta: \pi^*(TS^2,j) \longrightarrow (TE^v,\J)$ is
a $\C$-antilinear homomorphism.
Let $\U$ be the space of such $\theta$ whose $\floermetric$-norm is finite.

Consider first the condition that $D_{\jhat}(s)$ is onto. The main step
in the proof is to show that
\begin{align*}
\duniv(s,\jhat): \smooth(s^*TE^v) \times \U &\longrightarrow
\Omega^{0,1}(s^*TE^v,s^*\J),\\
\duniv(s,\jhat)(S,\theta) &= D_{\jhat}(s)S + \frac12 (\theta \circ D\pi)
\circ ds \circ j
\end{align*}
is onto for all $\jhat \in \modhat$ and $s \in \jjsec$.
$D\pi \circ ds = \id$ and therefore
$\duniv(s,\jhat)(0,\theta)_z = \frac12 \theta(s(z)) \circ j$. Because
$s: S^2 \longrightarrow E$ is an embedding, this means that
$\duniv(s,\jhat)(0 \times \U)$ is dense in $\Omega^{0,1}(s^*TE^v,s^*\J)$.
Since $D_{\jhat}(s)$ is a Fredholm operator, it follows that
$\duniv(s,\jhat)$ is onto.

Now consider the condition that $\ev$ is transverse to $\eta$. The
proof does not use any specific properties of $\eta$; for any
smooth map $c: C \longrightarrow E$, the space of $\jhat$ such that
$\ev$ is transverse to $c$ is dense in $\modhat$. As in
\cite[Section 6.1]{mcduff-salamon}, this follows from the fact
that the evaluation map on a suitable `universal moduli space'
is a submersion. The main step is to prove that the operator
\begin{align*}
T_zS^2 \times \smooth(s^*TE^v) \times \U &
\longrightarrow T_{s(z)}E \times \Omega^{0,1}(s^*TE^v,s^*\J),\\
(\xi, S, \theta) &\longmapsto (ds(z)\xi + S(z), \duniv(s,\jhat)(S,\theta))
\end{align*}
is onto for all $\jhat \in \modhat$, $s \in \jjsec$ and $z \in S^2$.
This is slightly stronger than the surjectivity of $\duniv$,
but it can be proved by the same argument. The proof of the
transversality condition for $\ev_{z_0}$ is similar. \qed

From now on, we assume that $(E,\O)$ is Hamiltonian. Let
$\tildeo \in \Omega^2(E)$ be a closed form with $\tildeo|E_z = \O_z$
for all $z$.

\begin{lemma} \label{th:thurston} Let $\hatj$ be an almost complex
structure on $E$ such that $\hatj \in \modhat$ for some $j,\J$.
There is a two-form $\sigma$ on $S^2$ such that $\tildeo +
\pi^*\sigma$ is a symplectic form on $E$ which tames $\hatj$. \qed
\end{lemma}

This is a well-known method for constructing symplectic forms,
due to Thurston \cite{thurston76}. Lemma \ref{th:thurston} makes
it possible to apply Gromov's compactness theorem to $(j,\hatj)$-holomorphic
sections, since any such section is a $\hatj$-holomorphic curve in $E$.
It is easy to see that any `bubble' in the limit of a sequence of 
$(j,\hatj)$-holomorphic sections must be a $J_z$-holomorphic curve in some
fibre $E_z$. As in \cite{ruan-tian94}, we can simplify the limits
by deleting some components and replacing multiply covered curves by
the underlying simple ones. If $\J \in \modeoreg$, every $J_z$ is
semi-positive, and the total Chern number does not increase
during the process. The outcome can be summarized as follows:

\begin{lemma} \label{th:bubbling} Let $j$ be a positively oriented
complex structure on $S^2$, $\hatj$ an almost complex structure
on $E$ such that $\hatj \in \modhat$ for some $\J \in \modeoreg$,
and $J$ the almost complex structure on $M$ which corresponds to $J_{z_0}$
through the isomorphism $i$.
Let $(s_m)_{m \in \N}$ be a sequence in $\jjsec$ such that
$\tildeo(s_m) \leq C$ and $c_1(TE^v,\O)(s_m) \leq c$ for all $m$.
Assume that $(s_m)$ has no convergent subsequence, and that $s_m(z_0)$
converges to $y \in E_{z_0}$. Then one of the following 
(not mutually exclusive) possibilities holds:
\begin{enumerate} \item \label{item:cancel-all-bubbles}
there is an $s \in \jjsec$ with $c_1(TE^v,\O)(s) < c$ and $s(z_0) = y$;
\item \label{item:bubble-one}
there are $s \in \jjsec$, $z \in S^2$ and 
$w \in \moduli^s_0(\J) \times_{PSL(2,\C)} \CP{1}$ such that
$c_1(TE^v,\O)(s) = c$, $s(z_0) = y$ and $s(z) = \eta(w)$;
\item \label{item:bubble-two} 
there is an $s \in \jjsec$ and a $v \in \cusp$ with
$k + c_1(TE^v,\O)(s) \leq c$, $i(\eta_1(v)) = s(z_0)$ and
$i(\eta_2(v)) = y$. 
\end{enumerate} \end{lemma}

\ref{item:bubble-two} is the case where bubbling occurs at $z_0$
in such a way that $y$ does not lie in the image of the principal
component of the limiting cusp-curve.

We will denote by $\sec(j,\hatj,S) \subset \jjsec$ the space of
$(j,\hatj)$-holomorphic sections of $E$ which lie in a given
$\Gamma$-equivalence class $S$.

\begin{lemma} \label{th:finite-non-empty}
For every $C \in \R$, there are only finitely
many $\Gamma$-equivalence classes $S$ with $\tildeo(S) \leq C$
and $\sec(j,\hatj,S) \neq \emptyset$. \end{lemma}

\proof The stronger statements in which $\Gamma$-equivalence
is replaced by homological equivalence or even homotopy are well-known (see 
\cite[Corollary 4.4.4]{mcduff-salamon}) consequences of 
Gromov's compactness theorem. \qed

\begin{proposition} \label{th:pseudo-cycle}
If $\hatj \in \modhatreg$ for some $j$ and $\J \in \modeoregz$, the map
\[
\ev_{z_0}: \sec(j,\hatj,S) \longrightarrow M, \quad
\ev_{z_0}(s) = i^{-1}(s(z_0))
\]
is a pseudo-cycle (in the sense of \cite[Section 7.1]{mcduff-salamon})
of dimension $d(S) = 2n + 2c_1(TE^v,\O)(S)$ for any $\Gamma$-equivalence class
$S$. \end{proposition}

\proof Because $\hatj$ is regular, $\sec(j,\hatj,S)$ is a smooth manifold.
Its dimension is given by the index of $D_{\jhat}(s)$ for any
$s \in \sec(j,\hatj,S)$, which is $d(S)$. It remains to show that $\ev_{z_0}$ 
can be compactified by countably many 
images of manifolds of dimension $\leq d(S) - 2$.
The compactification described in Lemma \ref{th:bubbling} has this
property: in case \ref{item:cancel-all-bubbles}, $y$ lies in the image of
\[
\ev_{z_0}: \sec(j,\hatj,S') \longrightarrow M
\]
for some $S'$ with $c_1(TE^v,\O)(S') < c_1(TE^v,\O)(S)$. 
Clearly $\dim \sec(j,\hatj,S') \leq d(S) -2$. Now, consider case
\ref{item:bubble-one}. Let $\sec_0$ be the space of $(j,\hatj)$-holomorphic
sections $s$ with $c_1(TE^v,\O)(s) = c_1(TE^v,\O)(S)$. The space
of $(s,z,w)$ as in Lemma \ref{th:bubbling}\ref{item:bubble-one}
is the inverse image of the diagonal $\Delta \subset E \times E$ by
the map
\[
\ev \times \eta: (S^2 \times \sec_0) \times
(\moduli^s_0(\J) \times_{PSL(2,\C)} \CP{1}) \longrightarrow E \times E.
\]
$\dim \moduli^s_0(\J) \times_{PSL(2,\C)} \CP{1} = \dim E - 4$,
and $S^2 \times \sec_0 $ is a manifold of dimension $d(S) + 2$.
Since $\hatj$ is regular, $\ev$ is transverse to $\eta$; it follows
that $(\ev \times \eta)^{-1}(\Delta)$ has dimension $d(S) - 2$.
Case \ref{item:bubble-two} is similar. \qed

The following way of assigning a homology class to a pseudo-cycle
is due to Schwarz \cite{schwarz96}.

\begin{lemma} \label{th:morse-homology}
Let $c: C \longrightarrow M$ be a $k$-dimensional
pseudo-cycle in $M$ which is compactified by
$c_\infty: C_\infty^{k-2} \longrightarrow M$. If $h$ is a Riemannian
metric on $M$ and $f \in \smooth(M,\R)$ is a Morse function for which
the Morse complex $(CM_*(f),\partial(f,h))$ is defined and
such that stable manifold $W^s(y;f,h)$ (for the negative gradient flow)
of any critical point $y$ of $f$ is transverse to $c,c_\infty$, the
sum
\[
\sum_y \#c^{-1}(W^s(y;f,h))\gen{y}
\]
over all critical points $y$ of Morse index $k$ is a cycle in $CM_*(f)$
(we use Morse complexes with $\Z/2$-coefficients). \qed \end{lemma}

Such a pair $(f,h)$ always exists, and the homology class of the
cycle in the `Morse homology' \cite{schwarz} of $M$ is independent
of all choices. Since Morse homology is canonically isomorphic to
singular homology \cite{schwarz97}, this defines a class in
$H_k(M;\Z/2)$ which we denote by $[c(C)]$. If $C$ is compact,
$[c(C)]$ is the fundamental class in the classical sense.

\begin{defn} \label{def:bigq}
Let $(E,\O,S)$ be a normalized Hamiltonian fibre bundle over 
$S^2$ with fibre $\mo$. Let $d = 2n + 2c_1(TE^v,\O)(S)$. We define
\[
Q(E,\O,S) = \sum_{\gamma \in \Gamma}
[\ev_{z_0}(\sec(j,\hatj,\gamma + S))] \otimes \gen{\gamma}
\; \in QH_d\mo,
\]
where $j$ is a positively oriented complex structure on $S^2$
and $\hatj \in \modhatreg$ for some $\J \in \modeoregz$. \end{defn}

The notation $\gamma + S$ was introduced in section \ref{sec:elementary}.
The formal sum defining $Q(E,\O,S)$ lies in $\quantum$ because of Lemma
\ref{th:finite-non-empty}. It is of degree $d$ because the dimension
of $\sec(j,\hatj,\gamma + S)$ is $d + 2c_1(\gamma)$ and $\gen{\gamma}$
has degree $-2c_1(\gamma)$.

\begin{prop} $Q(E,\O,S)$ is independent of the choice of $j$,
$\J$ and $\hat{J}$. \end{prop}

We omit the proof. It is based on the fact that cobordant pseudo-cycles
determine the same homology class, and is similar to
\cite[Proposition 7.2.1]{mcduff-salamon}.

Sometimes it is convenient to define $QH_*\mo$ in terms of
Morse homology as the homology of the graded tensor
product $(CM_*(f) \otimes \Lambda, \partial(f,h) \otimes \id)$.
An element of $CM_*(f) \otimes \Lambda$ is a (possibly infinite)
linear combination of $\gen{y} \otimes \gen{\gamma}$ for
$y \in \Crit(f)$, $\gamma \in \Gamma$.
Assume that $(f,h)$ satisfies the conditions of Lemma \ref{th:morse-homology}
with respect to the pseudo-cycles 
$ev_{z_0}: \sec(j,\hatj,S') \longrightarrow M$ for all $S'$, and define
\begin{equation*} 
\sec(j,\hatj,S',y) = \{ s \in \sec(j,\hatj,S') \; | \;
\ev_{z_0}(s) \in W^s(y;f,h) \}
\end{equation*}
for $y \in \Crit(f)$. Then $Q(E,\O,S)$ is the homology class
of the cycle
\begin{equation} \label{eq:explicit}
\sum_{y,\gamma} \#\sec(j,\hatj,\gamma + S,y) \gen{y} \otimes
\gen{\gamma} \in CM_*(f) \otimes \Lambda,
\end{equation}
where the sum is over all $(y,\gamma) \in \Crit(f) \times \Gamma$
such that $i_f(y) = d + 2c_1(\gamma)$ ($i_f$ denotes the Morse index).

It is often difficult to decide whether a given $\hatj$ and its
restriction $\J = \hatj|TE^v$ are regular, since this depends on
all holomorphic sections and the holomorphic curves in the fibres. Moreover, 
in many examples the most natural choice of $\hatj$ is not regular.
However, almost complex structures which satisfy a weaker condition
can be used to determine $Q(E,\O,S)$ partially: to
compute $[\ev_{z_0}(\sec(j,\hatj,\gamma + S))]$ for a single $\gamma$, 
it is sufficient that $\sec(j,\hatj,\gamma + S)$ and all its possible limits
listed in Lemma \ref{th:bubbling} are regular. 
Since the contributions from those $\gamma$ with $d + 2c_1(\gamma) < 0$ 
or $d + 2c_1(\gamma) > 2n$ are zero, sometimes all of 
$Q(E,\O,S)$ can be computed using an almost complex structure 
which is not regular, as in the following case:

\begin{proposition} \label{th:computation}
Let $(E,J)$ be a compact complex manifold, $\pi: E
\longrightarrow \CP{1}$ a holomorphic map with no critical points,
$\tildeo \in \Omega^2(E)$ a closed form whose restrictions
$\O_z = \tildeo|E_z$ are K{\"a}hler forms,
and $i: \! \mo \longrightarrow (E_{z_0},\O_{z_0})$ an
isomorphism for some $z_0 \in \CP{1}$.  Assume that
\begin{enumerate} 
\item the space $\sec$ of holomorphic sections $s$ of $\pi$ 
with $c_1(TE^v)(s) \leq 0$ is connected. In particular, all
of them lie in a single $\Gamma$-equivalence class $S_0$.
\item For any $s \in \sec$, $H^{0,1}(\CP{1},s^*TE^v) = 0$.
\item Any holomorphic map $w \in \smooth(\CP{1},E)$ such that
$\im(w) \subset E_z$ for some $z \in \CP{1}$ satisfies $c_1(TE)(w) \geq 0$;
\item if $w$ is as before and not constant, and
$c_1(TE)(w) + c_1(TE^v)(S_0) \leq 0$,
then $s(z) \notin \im(w)$ for all $s \in \sec$.
\end{enumerate}
In that case, $\sec$ is a smooth compact manifold and
\[
Q(E,\O,S_0) = (\ev_{z_0})_*[\sec] \otimes \gen{0},
\]
where $\gen{0} \in \novi$ is the element corresponding to the
trivial class in $\Gamma$.
\end{proposition}

\newcommand{\pplus}{\Psi^+(e)}
\newcommand{\pairsofmaps}{\mathcal{P}(\gamma)}
\newcommand{\choices}{\mathcal{C}(H^{\infty},\J^{\infty},
H^{-\infty},\J^{-\infty})}
\newcommand{\pairs}{\mathcal{S}'(H^+,\J^+,H^-,\J^-,\gamma,y)}
\newcommand{\approxiglue}{\widehat{\#}_R}
\newcommand{\norm}{\|}

\newcommand{\noduli}{?? \mathcal{N}}
\section{A gluing argument \label{sec:gluing}}

\begin{definition} For $\pair \in \group$, define
\[
q\pair = Q(E_g,\O_g,S_{\tilde{g}}) \in \quantum.
\]
\end{definition}

By definition, $Q(E_g,\O_g,S_{\tilde{g}})$ is an element
of degree $2n + 2c_1(TE_g^v,\O_g)(S_{\tilde{g}})$. Using
\eqref{eq:maslov-and-chern} it follows that
\[
q\pair \in QH_{2n - 2I\pair}\mo.
\]

Let $e = [M] \otimes \gen{0}$ be the `fundamental class' in
$\quantum$. The connection between the
invariant $q$ and the $\group$-action on Floer homology is given by

\begin{theorem} \label{th:gluing}
For all $(g,\tilde{g}) \in \group$,
\[
q\pair = \Psi^- \, HF_*(g,\tilde{g}) \, \pplus.
\] \end{theorem}

Here
\[
\Psi^+: QH_*\mo \longrightarrow HF_*\mo, \quad 
\Psi^-: HF_*\mo \longrightarrow QH_*\mo
\]
are the canonical homomorphisms of Piunikhin, Salamon and Schwarz {\pss}. 
They proved that $\Psi^+,\Psi^-$ are isomorphisms and inverses of each 
other, and a part of their proof will serve as a model for the 
proof of Theorem \ref{th:gluing}. More precisely, we consider 
the proof of $\id_{QH_*\mo} = \Psi^- \Psi^+$ given in {\pss} 
and specialize it to the case
\[
e = \Psi^- \pplus.
\]
There are two steps in proving this: first, a gluing argument shows 
that $\Psi^- \pplus$ is equal to a certain Gromov-Witten invariant.
Then one computes that the value of this invariant is again $e$. 
The Gromov-Witten invariant which arises can be described in our
terms as the $Q$-invariant of the trivial Hamiltonian fibre
bundle $(S^2 \times M,\o)$ with the $\Gamma$-equivalence class $S_0$
containing the constant sections. Therefore, what the gluing 
argument in {\pss} proves is
\begin{equation} \label{eq:pss}
Q(S^2 \times M, \omega, S_0) = \Psi^- \pplus.
\end{equation}
This is precisely the special case $(g,\tilde{g}) = \id_{\group}$ of
Theorem \ref{th:gluing}, since in that case $(E_g,\O_g,S_{\tilde{g}}) = 
(S^2 \times M,\o,S_0)$ and $HF_*(g,\tilde{g})$ 
is the identity map. The analytical details
of the proof of Theorem \ref{th:gluing} are practically
the same as in the special case \eqref{eq:pss}; the
only difference is that the gluing procedure is modified
by `twisting' with $g$. Since a detailed account of the results
of {\pss} based on the analysis in \cite{schwarz95} is in
preparation, we will only sketch the proof of Theorem \ref{th:gluing}.

The basic construction underlying the proof is familiar from
section \ref{sec:elementary}: a pair $(v^+,v^-) \in
\smooth(D^+,M) \times \smooth(D^-,M)$ such that $v^+|\partial D^+ = x
$ and $v^-|\partial D^- = g(x)$ for some $x \in \loops$ 
gives rise to a section of $E_g$. If
$[v^-,g(x)] = \tilde{g}([v^+,x]) \in \tloops$, 
this section lies in the $\Gamma$-equivalence class
$S_{\tilde{g}}$. More generally, if $[v^-,g(x)] =
\gamma \cdot \tilde{g}([v^+,x])$ for some $\gamma \in \Gamma$,
the section lies in $(-\gamma) + S_{\tilde{g}}$. The sign
$-\gamma$ occurs for the following reason: to obtain the element
$[v^-,g(x)] \in \tloops$, one uses a diffeomorphism 
from $D^-$ to the standard disc $D^2$. Following Convention 
\ref{th:orientation-convention}, this diffeomorphism should be 
orientation-reversing.

To adapt this construction to the framework of surfaces with
tubular ends, we replace $D^+,D^-$ by $\Sigma^+ = D^+ \cup_{S^1}
(\R^+ \times S^1)$ and $\Sigma^- = (\R^- \times S^1) \cup_{S^1}
D^-$, with the orientations induced from the standard ones on
$\R^{\pm} \times S^1$. A map $u \in \smooth(\Sigma^+,M)$ such that
$\lim_{s \rightarrow \infty} u(s,\cdot) = x$ for some $x \in \loops$
can be extended continuously to the compactification $\overline{\Sigma}^+ =
\Sigma^+ \cup (\{\infty\} \times S^1)$. Using an orientation-preserving
diffeomorphism $D^2 \longrightarrow \overline{\Sigma}^+$, this defines
an element $c = [u,x] \in \tloops$. We say that $u$ converges to $c$.
There is a parallel notion for $u \in \smooth(\Sigma^-,M)$, except
that in this case the identification of $\overline{\Sigma}^-$ with $D^2$
reverses the orientation.

\begin{defn} For $\gamma \in \Gamma$, let $\pairsofmaps$ be the space of
pairs $(u^+,u^-) \in \smooth(\Sigma^+,M) \times \smooth(\Sigma^-,M)$ such
that $u^+$ converges to $c$ and $u^-$ to $(-\gamma) \cdot \tilde{g}(c)$
for some $c \in \tloops$.
\end{defn}

Fix a point $z_0 \in D^- \subset \Sigma^-$. We will denote the trivial
symplectic fibre bundles $\Sigma^\pm \times \mo$ by $(E^\pm,\O^\pm)$,
and the obvious isomorphism $\mo \longrightarrow (E^-_{z_0},\O^-_{z_0})$
by $i$. For $R \geq 0$, consider $\Sigma_R^+ = D^+ \cup_{S^1}
([0;R] \times S^1) \subset \Sigma^+$, $\Sigma_R^- = 
([-R;0] \times S^1) \cup D^- \subset \Sigma^-$ and $\Sigma_R =
\Sigma_R^+ \cup_{\partial\Sigma_R^{\pm}} \Sigma_R^-$, with the
induced orientation. With respect to Riemannian metrics
on $\Sigma^+,\Sigma^-$ whose restriction to $\R^{\pm} \times S^1$
is the standard tubular metric, $\Sigma_R$ is a surface with
an increasingly long `neck' as $R \rightarrow \infty$. Let
$(E_R,\O_R)$ be the symplectic fibre bundle over $\Sigma_R$ 
obtained by gluing together $E^+_R = E^+|\Sigma^+_R$ and
$E^-_R = E^-|\Sigma^-_R$ using $\phi_g:
E^+|\partial\Sigma^+_R \longrightarrow E^-|\partial\Sigma_R^-$, 
$\phi_g(R,t,y) = (-R,t,g_t(y))$. The terminology introduced in section 
\ref{sec:elementary} applies to bundles over $\Sigma_R$ since it is
homeomorphic to $S^2$. $(E_R,\O_R)$ is Hamiltonian; 
it differs from $(E_g,\O_g)$ only because the base is parametrized 
in a different way (more precisely, $(E_R,\O_R)$ is the pullback of 
$(E_g,\O_g)$ by an oriented diffeomorphism
$\Sigma_R \longrightarrow S^2$). Let $S_R$ be the $\Gamma$-equivalence
class of sections of $E_R$ which corresponds to $S_{\tilde{g}}$.

For $(u^+,u^-) \in \pairsofmaps$, one can construct a section
$u^+ \approxiglue u^-$ of $E_R$ for large $R$ by
`approximate gluing': let $x = \lim_{s \rightarrow \infty} u^+(s,\cdot) 
\in \loops$. There is an $R_0 \geq 0$ and a family 
$(\hat{u}_R^+,\hat{u}_R^-)_{R \geq R_0}$ in $\pairsofmaps$ 
which converges uniformly to $(u^+,u^-)$ as $R \rightarrow \infty$
and such that $\hat{u}_R^+(s,t) = x(t)$ for $s \geq R$,
$\hat{u}_R^-(s,t) = g_t(x(t))$ for $s \leq -R$. Define
\[
(u^+ \approxiglue u^-)(z) = \begin{cases}
(z,\hat{u}_R^+(z)) & z \in \Sigma_R^+,\\
(z,\hat{u}_R^-(z)) & z \in \Sigma_R^-.
\end{cases}
\]
If $R$ is large, $u^+ \approxiglue u^-$ lies in the $\Gamma$-equivalence 
class $\gamma + S_R$.

We will now introduce nonlinear $\bar{\partial}$-equations 
on $\Sigma^+,\Sigma^-$ and state the gluing theorem for solutions
of them. Like the pair-of-pants product, these equations are part
of the formalism of `relative Donaldson type invariants' of {\pss}.

Let $j^+$ be a complex structure on $\Sigma^+$ whose restriction
to $\R^+ \times S^1$ is the standard complex structure. Choose
a regular pair $(H^{\infty},\J^{\infty})$ and
$H^+ \in \smooth(\R^+ \times S^1 \times M,\R),
\J^+ \in \modJ{\Sigma^+}$ such that $H^+|[0;1] \times S^1 \times M = 0$
and $H^+(s,t,\cdot) = H^{\infty}(t,\cdot), J^+_{(s,t)} = J^{\infty}_t$ for 
$s \geq 2$, $t \in S^1$. For $c \in \Crit(a_{H^{\infty}}) \subset \tloops$, 
let $\moduli^+(c;H^+,\J^+)$ be the space of $u \in \smooth(\Sigma^+,M)$
which converge to $c$ and satisfy
\begin{equation} \label{eq:plus} \begin{aligned}
du(z) + J_z^+ \circ du(z) \circ j^+ &= 0 \quad \text{ for } z \in D^+,\\
\partials{u} + J_z^+(u)\left(\partialt{u} - X_{H^+}(s,t,u)\right) &= 0
\quad \text{ for } z = (s,t) \in \R^+ \times S^1.
\end{aligned} \end{equation}
For generic $(H^+,\J^+)$, $\moduli^+(c;H^+,\J^+)$ is a manifold
of dimension $2n - \mu_{H^{\infty}}(c)$ and the zero-dimensional
spaces are finite.

To write down the corresponding equations for $u \in \smooth(\Sigma^-,M)$,
choose a complex structure $j^-$ on $\Sigma^-$,
a regular pair $(H^{-\infty},\J^{-\infty})$ and
$H^- \in \smooth(\R^- \times S^1 \times M,\R)$, 
$\J^- \in \modJ{\Sigma^-}$ with properties symmetric to those above. 
The equations are 
\begin{equation} \label{eq:minus} \begin{aligned}
du(z) + J_z^- \circ du(z) \circ j^- &= 0 \quad \text{ for } z \in D^-,\\
\partials{u} + J_z^-(u)\left(\partialt{u} - X_{H^-}(s,t,u)\right) &= 0
\quad \text{ for } z = (s,t) \in \R^- \times S^1.
\end{aligned} \end{equation}
Let $h$ be a Riemannian metric on $M$ and $f \in \smooth(M,\R)$ a 
Morse function. If $c$ is a critical point of $a_{H^{-\infty}}$ and $y$ 
a critical point of $f$, we denote by $\moduli^-(c,y;H^-,\J^-)$ the space of
solutions $u$ of \eqref{eq:minus} which converge to $c$ and with
$u(z_0) \in W^s(y;f,h)$. In the generic case, $\moduli^-(c,y;H^-,\J^-)$
is a manifold of dimension $\mu_{H^{-\infty}}(c) - i_f(y)$, and the
zero-dimensional spaces are again finite. For fixed 
$(H^\infty,\J^\infty)$ and $(H^{-\infty},\J^{-\infty})$, we 
will denote the space of all $(H^+,\J^+,H^-,\J^-,f,h)$ by
$\choices$.

\eqref{eq:plus} and \eqref{eq:minus} can be written in a different
way using an idea of Gromov. For $(z,y) \in \Sigma^+ \times M$, let
$\nu^+(z,y): T_z\Sigma^+ \longrightarrow T_yM$ be the 
$(j^+,J_z^+)$-antilinear homomorphism given by
\[
\nu^+(z,y) = \begin{cases}
ds \otimes J_z^+ X_{H^+}(s,t,y) + dt \otimes X_{H^+}(s,t,y) &\text{for }
z = (s,t),\\
0 &\text{for } z \in D^+. \end{cases}
\]
$u$ is a solution of \eqref{eq:plus} iff
\begin{equation} \label{eq:plus-two}
du(z) + J_z^+ \circ du(z) \circ j^+ = \nu^+(z,u(z)).
\end{equation}
One can think of $\J^+ = (J^+_z)_{z \in \Sigma^+}$ as a family of almost
complex structures on the fibres of the trivial bundle $E^+$. 
Consider the almost complex structure
\[
\hatj^+_{(z,y)}(Z,Y) = (j^+Z,J_z^+Y + \nu^+(z,y)(j^+Z))
\]
on $E^+$. $\hatj^+| \{z\} \times M = J^+_z$ for all $z \in \Sigma^+$, and 
the projection $E^+ \longrightarrow \Sigma^+$ is $(\hatj^+,j^+)$-linear. 
A straightforward computation shows that $u$ is a solution of 
\eqref{eq:plus-two} iff the section $s(z) = (z,u(z))$ of $E^+$ is 
$(j^+,\hatj^+)$-holomorphic. There is an almost
complex structure $\hatj^-$ on $E^-$ such that solutions
of \eqref{eq:minus} correspond to $(j^-,\hatj^-)$-holomorphic sections
in the same way.

From now on, we will assume that $(H^{\infty},\J^{\infty})$ is the pullback
of $(H^{-\infty},\J^{-\infty})$ by $g$. As a consequence, 
$(J^+_z)_{z \in \Sigma^+_R}$ and $(J^-_z)_{z \in \Sigma^-_R}$ can
be pieced together to a smooth family $\J_R = (J_{R,z})_{z \in \Sigma_R}
\in \mathcal{J}(E_R,\O_R)$ for any $R \geq 2$. $j^+|\Sigma^+_R$ and
$j^-|\Sigma^-_R$ determine a complex structure $j_R$ on
$\Sigma_R$, and there is a unique $\hatj_R \in 
\widehat{\mathcal{J}}(j_R,\J_R)$ such that $\hatj_R|E^{\pm}_R =
\hatj^{\pm}|E^{\pm}_R$. As in the previous section, we denote
the space of $(j_R,\hatj_R)$-holomorphic sections of $E_R$ in a 
$\Gamma$-equivalence class $S$ by $\sec(j_R,\hatj_R,S)$, and by
$\sec(j_R,\hatj_R,S,y)$ for $y \in \Crit(f)$ the subspace
of those $s$ such that $i^{-1}(s(z_0)) \in W^s(y;f,h)$.

For $(\gamma,y) \in \Gamma \times \Crit(f)$ with
\begin{equation} \label{eq:dimension-condition}
i_f(y) = 2n - 2I\pair + 2c_1(\gamma),
\end{equation}
let $\pairs$ be the disjoint union of 
\begin{equation} \label{eq:the-product}
\moduli^+(c;H^+,\J^+) \times \moduli^-( (-\gamma) \cdot \tilde{g}(c),y;
H^-,\J^-)
\end{equation}
for all $c \in \Crit(a_{H^\infty})$. Equivalently, $\pairs$ is the
space of $(u^+,u^-) \in \pairsofmaps$ such that $s^+(z) = (z,u^+(z))$ is 
$(j^+,\hatj^+)$-holomorphic, $s^-(z) = (z,u^-(z))$ is 
$(j^-,\hatj^-)$-holomorphic and $i^{-1}(s^-(z_0)) \in W^s(y;f,h)$.

For generic $(H^+,\J^+)$ and $(H^-,\J^-)$,
$\dim \moduli^+(c;H^+,\J^+) = 2n - \mu_{H^\infty}(c)$ and
$
\dim \moduli^-((-\gamma) \cdot \tilde{g}(c),y;H^-,\J^-) 
= \mu_{H^{-\infty}}((-\gamma) \cdot \tilde{g}(c)) - i_f(y)
= \mu_{H^\infty}(c) - 2n
$
(the last equality uses \eqref{eq:dimension-condition}, 
\eqref{eq:grading-shift} and Lemma \ref{th:shift}). Therefore the product
\eqref{eq:the-product} is zero-dimensional 
if $\mu_{H^\infty}(c) = 2n$, and is empty otherwise. 
A compactness theorem shows that $\pairs$ is a finite set.

We can now state the `gluing theorem' which is the main
step in the proof of Theorem \ref{th:gluing}.

\begin{theorem} \label{th:technical-gluing-theorem}
There is a subset of $\choices$ of second category such that
if $(H^+,\J^+,H^-,\J^-,f,h)$ lies in this subset, the following
holds: for $(\gamma,y) \in \Gamma \times \Crit(f)$ satisfying
\eqref{eq:dimension-condition}, there is an $R_0 > 2$ and a family
of bijective maps
\[
\glue_R: \pairs \longrightarrow \sec(j_R,\hatj_R,\gamma + S_R,y)
\]
for $R \geq R_0$. \end{theorem}

The construction of $\glue_R$ relies on the following property
of the `approximate gluing' $\approxiglue$ for $(u^+,u^-) \in \pairs$:
because $\hat{u}_R^+$ converges uniformly to $u^+$ as
$R \rightarrow \infty$, it approximately satisfies
\eqref{eq:plus} for large $R$. Therefore 
$\hat{s}^+_R(z) = (z,\hat{u}_R^+(z))$
is an approximately $(j^+,\hatj^+)$-holomorphic section of $E^+$.
Using the same argument for $u^-$, it follows that
$\hat{s}_R = u^+ \approxiglue u^-$ is an approximately
$(j_R,\hatj_R)$-holomorphic section of $E_R$. More precisely,
\[
\norm d\hat{s}_R + \hatj_R \circ d\hat{s}_R \circ 
j_R \norm_p \rightarrow 0 \quad \text{ as } R \rightarrow \infty.
\]
Here $\norm \cdot \norm_p$ $(p>2)$ is the $L^p$-norm on $\Sigma_R$ with 
respect to a family of metrics with an 
increasingly long `neck' as described above. 

$\approxiglue$ can be set up in such a way that $i^{-1}(\hat{s}_R(z_0))
\in W^s(y;f,h)$ for all $R$. The section $u^+ \glue_R u^-$
is obtained from $\hat{s}_R$ by an application of the implicit function
theorem in the space of $W^{1,p}$-sections $s$ of $E_R$ which satisfy
$i^{-1}(s(z_0)) \in W^s(y;f,h)$. The argument can be modelled on
\cite[Section 4.4]{schwarz95}. For large $R$,
$u^+ \glue_R u^-$ and $\hat{s}_R$ are close in the $W^{1,p}$-metric
and therefore homotopic. It follows that $u^+ \glue_R u^-$ lies
in the $\Gamma$-equivalence class $\gamma + S_R$.

By construction, $\glue_R$ is injective for large $R$.
Surjectivity is proved by a limiting argument (`stretching the neck').
Let $s_m \in \sec(j_{R_m},\hatj_{R_m},\gamma + S_{R_m},y)$ be a
sequence of sections, with $R_m \rightarrow \infty$, such that
$s_m | \Sigma^+_{R_m}$ and $s_m | \Sigma^-_{R_m}$ converge on compact subsets
to sections $s^+(z) = (z,u^+(z))$ of $E^+$ and $s^-(z) = (z,u^-(z))$
of $E^-$. Then $s^\pm$ is $(j^\pm,\hatj^\pm)$-holomorphic.
If the pair $(s^+,s^-)$ describes the `geometric limit' of the 
sequence $s_m$ completely, $(u^+,u^-) \in \pairs$ and
$s_m = u^+ \glue_{R_m} u^-$ for large $m$. Transversality and
dimension arguments are used to exclude more complicated limiting
behaviour.

It remains to explain how Theorem \ref{th:gluing} is derived
from the technical Theorem \ref{th:technical-gluing-theorem}.
The main point is that $\Psi^+,\Psi^-$ are the `relative
Donaldson type invariants' associated to $\Sigma^+,\Sigma^-$.
The element $\pplus \in HF_{2n}\mo$ has a particularly simple
description in terms of the moduli spaces $\moduli^+$: it is
the homology class of the cycle
\[
\Psi^+(e;H^+,\J^+) = \sum_{c \in \Crit_{2n}(a_{H^\infty})} 
\#\moduli^+(c;H^+,\J^+) \gen{c} \in CF_{2n}(H^\infty)
\]
for generic $(H^+,\J^+)$. It follows that
$HF_*\pair \pplus \in HF_{2n - 2I\pair}\mo$ is represented by
\[
\sum_c \#\moduli^+(c;H^+,\J^+)\gen{\tilde{g}(c)} 
\in CF_{2n - 2I\pair}(H^{-\infty}).
\]
To define $\Psi^-$, we consider $\quantum$ as the homology
of $(CM_*(f) \otimes \novi, \partial(f,h) \otimes \id)$.
For generic $(H^-,\J^-)$, the formula
\[
\Psi^-(H^-,\J^-,f,h)(\gen{c}) = \sum_{\gamma,y}
\#\moduli^-(\gamma \cdot c,y;H^-,\J^-) \gen{y} \otimes \gen{\!-\!\gamma},
\]
where the sum is over all $(\gamma,y) \in \Gamma \times \Crit(f)$ such that
$i_f(y) = \mu_{H^{-\infty}}(c) - 2c_1(\gamma)$, defines a $\novi$-linear 
homomorphism of chain complexes 
\[
\Psi^-(H^-,\J^-,f,h): CF_*(H^{-\infty}) \longrightarrow 
CM_*(f) \otimes \novi. 
\]
$\Psi^-$ is defined as the induced map on homology groups. 
Here, as in the case of $\Psi^+$, we have omitted
several steps in the construction; we refer to {\pss} for a more
complete description. As usual, our setup differs slightly from
that in {\pss} because the almost complex structures may be $z$-dependent; 
however, the proofs can be easily adapted, provided always that 
{\wplus} is satisfied.

By comparing the expressions for $HF_*(g,\tilde{g}) \pplus$ and
$\Psi^-$ with the definition of $\pairs$, one obtains
\[
\Psi^- HF_*\pair \pplus = 
\sum_{\gamma \in \Gamma} [c_\gamma] \otimes \gen{\gamma},
\]
where $[c_\gamma] \in H_*(M;\Z/2)$ is the homology class of the
cycle
\[
c_\gamma = \sum_y \#\pairs \gen{y} \in CM_*(f);
\]
the sum runs over all $y \in \Crit(f)$ with
$i_f(y) = 2n - 2I\pair + 2c_1(\gamma)$.

Clearly, $Q(E_R,\O_R,S_R) = Q(E_g,\O_g,S_{\tilde{g}})$ for all $R \geq 0$.
To obtain an explicit expression for $Q(E_R,\O_R,S_R)$, choose 
$\J_R' \in \mathcal{J}^{\mathrm{reg},z_0}(E_R,\O_R)$,
$\hatj_R' \in \widehat{\mathcal{J}}^{\mathrm{reg}}(j_R,\J_R')$, a Riemannian
metric $h'$ and a Morse function $f'$, all of which satisfy 
appropriate regularity conditions. By \eqref{eq:explicit},
\[
Q(E_R,\O_R,S_R) = \sum_{\gamma \in \Gamma} [c_{\gamma,R}] 
\otimes \gen{\gamma},
\]
where $c_{\gamma,R} \in CM_*(f)$ is given by
\[
c_{\gamma,R} = \sum_{y'} \#\sec(j_R,\hatj_R',\gamma + S_R,y') \gen{y'}.
\]
This time the sum is over those $y' \in \Crit(f')$ such that
$i_{f'}(y') = 2n + 2c_1(TE_R^v,\O_R)(S_R) + 2c_1(\gamma)$, but
by equation \eqref{eq:maslov-and-chern} the Morse indices
are the same as above.

The statement of Theorem \ref{th:gluing} is that
$[c_{\gamma,R}] = [c_\gamma]$ for all $\gamma$.
Since $[c_{\gamma,R}]$ is independent of $R$, it is sufficient
to show that for each $\gamma$ there is an $R$ such that 
for a suitable choice of $\J_R',\hatj_R',f',h'$,
$[c_{\gamma,R}] = [c_\gamma]$. Choose some $\gamma_0 \in \Gamma$. 
We can assume that $H^+,\J^+,H^-,\J^-,f,h$ have been chosen 
as in Theorem \ref{th:technical-gluing-theorem}. Since $f$ has only
finitely many critical points, there is an $R$ such that
Theorem \ref{th:technical-gluing-theorem} holds for all $(\gamma_0,y)$.
Then
\[
c_{\gamma_0} = \sum_y \#\sec(j_R,\hatj_R,\gamma_0 + S_R,y) \gen{y}.
\]
and it follows that $c_{\gamma_0} = c_{\gamma_0,R}$ if we take
$\J_R' = \J_R$, $\hatj_R' = \hatj_R$ and $(f',h') = (f,h)$. 
Note that $\J_R \in \mathcal{J}(E_R,\O_R)$ and
$\hatj_R \in \widehat{\mathcal{J}}(j_R,\J_R)$ 
are not `generic' choices because they are $s$-independent on the
`neck' of $\Sigma_R$. However, a final consideration shows that for
large $R$, they can be used to compute the coefficient 
$[c_{\gamma_0,R}]$ of $Q(E_R,\O_R,S_R)$.

\newcommand{\gw}{\tilde{\Phi}_{(A,\o)}}
\section{Proof of the main results \label{sec:proofs}} 

For $a_1, a_2, a_3 \in H_*(M;\Z/2)$ and $A \in H_2(M;\Z)$, let
$\gw(a_1,a_2,a_3) \in \Z/2$ be the mod $2$ reduction of the
Gromov-Witten invariant of \cite[Section 8]{ruan-tian94}.
It is zero unless
\begin{equation} \label{eq:dimension-of-gw}
\dim(a_1) + \dim(a_2) + \dim(a_3) + 2c_1(A) = 4n.
\end{equation}
Intuitively, $\gw(a_1,a_2,a_3)$ is the number (modulo $2$) of
$J$-holomorphic spheres in the class $A$ which meet suitable 
cycles representing $a_1,a_2$ and $a_3$. 
Let $a_1 \qp_A a_2 \in H_*(M;\Z/2)$ be
the class defined by
\[
(a_1 \qp_A a_2) \cdot a_3 = \gw(a_1,a_2,a_3) \text{ for all } a_3,
\]
where $\cdot$ is the ordinary intersection product.
For $\gamma \in \Gamma$, let $a_1 \qp_{\gamma} a_2 \in H_*(M;\Z/2)$
be the sum of $a_1 \qp_A a_2$ over all classes $A$ which can
be represented by a smooth map $w: S^2 \longrightarrow M$ with
$[w] = \gamma$ (only finitely many terms of this sum are nonzero). 
The quantum intersection product $\qp$ on $\quantum$ is defined 
by the formula
\[
(a_1 \otimes \gen{\gamma_1}) \qp (a_2 \otimes \gen{\gamma_2}) =
\sum_{\gamma \in \Gamma} 
(a_1 \qp_\gamma a_2) \otimes \gen{\gamma_1 + \gamma_2 + \gamma},
\]
extended to infinite linear combinations in the obvious way. 
$\qp$ is a bilinear $\novi$-module map with the following properties:
\begin{enumerate}
\item \label{item:grad} if $a_1 \in QH_i\mo$ and $a_2 \in QH_j\mo$,
$a_1 \qp a_2 \in QH_{i+j-2n}\mo$. 
\item \label{item:as} $\qp$ is associative.
\item \label{item:un}
$e = [M] \otimes \gen{0} \in QH_{2n}\mo$ is the unit of $\qp$.
\item \label{item:pp}
$\Psi^+(a_1 \qp a_2) = \Psi^+(a_1) \pairofpants \Psi^+(a_2) \in HF_*\mo.$ 
\end{enumerate} 
\ref{item:grad} follows from \eqref{eq:dimension-of-gw}.
\ref{item:as} is due to Ruan and Tian \cite{ruan-tian94};
another proof is given in \cite{mcduff-salamon}. 
\ref{item:un} is \cite[Proposition 8.1.4(iii)]{mcduff-salamon}
and \ref{item:pp} is the main result of {\pss}. From the two last
items and the fact that $\Psi^+$ is an isomorphism, one obtains

\begin{lemma} \label{th:pp-unit}
$u = \Psi^+(e)$ is the unit of $(\floer,\pairofpants)$. \qed
\end{lemma}

\proof[Proof of Theorem \ref{bigth:comparison}] By Theorem \ref{th:gluing},
$q(g,\tilde{g}) = \Psi^- HF_*(g,\tilde{g})(u)$. Because $\Psi^-$ is the 
inverse of $\Psi^+$, this can be written as $\Psi^+(q(g,\tilde{g})) = 
HF_*(g,\tilde{g})(u)$. From Lemma \ref{th:pp-unit} and Proposition 
\ref{th:pair-of-pants-linear}, it follows that
\[ \begin{split}
HF_*\pair(b) &= HF_*\pair(u \pairofpants b)\\
&= HF_*\pair(u) \pairofpants b = \Psi^+(q\pair) \pairofpants b. \qed
\end{split} \]

\proof[Proof of Proposition \ref{bigth:trivial}] Proposition 
\ref{prop:further-properties}\ref{item:lambda-multiplication} 
says that $HF_*(\id,\gamma)$ is given by multiplication with 
$\gen{\gamma} \in \novi$. Using Theorem \ref{th:gluing}, 
$\Psi^-\Psi^+ = \id$ and the fact that 
$\Psi^-$ is $\Lambda$-linear, we can compute
\begin{multline*}
q(\id,\gamma) = \Psi^- HF_*(\id,\gamma)(u) = \Psi^- (\gen{\gamma} \, u) = \\
= \gen{\gamma} \Psi^-(u) = \gen{\gamma} ([M] \otimes \gen{0}) 
= [M] \otimes \gen{\gamma}. \qed
\end{multline*}

\proof[Proof of Corollary \ref{bigth:homomorphism}] Essentially, this
is a consequence of the fact that the maps $HF_*\pair$ define a
$\group$-action on Floer homology (Proposition 
\ref{prop:further-properties}\ref{item:functor}). Take
$(g_1,\tilde{g}_1), (g_2,\tilde{g}_2) \in \group$. Using
Theorem \ref{th:gluing} twice, we obtain
\[ \begin{split}
q(g_1g_2,\tilde{g}_1\tilde{g}_2) &= 
\Psi^- HF_*(g_1g_2,\tilde{g}_1\tilde{g}_2)(u) \\
& = \Psi^- HF_*(g_1,\tilde{g}_1) HF_*(g_2,\tilde{g}_2)(u)\\
& = \Psi^- HF_*(g_1,\tilde{g}_1) \Psi^+(q(g_2,\tilde{g}_2)).
\end{split} \]
Theorem \ref{bigth:comparison} with 
$(g,\tilde{g}) = (g_1,\tilde{g}_1)$ and $b =
\Psi^+(q(g_2,\tilde{g}_2))$ says that
\[
HF_*(g_1,\tilde{g}_1) \Psi^+(q(g_2,\tilde{g}_2)) = 
\Psi^+(q(g_1,\tilde{g}_1)) \pairofpants \Psi^+(q(g_2,\tilde{g_2})).
\]
Using property \ref{item:pp} above, we conclude that
\[
q(g_1g_2,\tilde{g}_1\tilde{g}_2) = q(g_1,\tilde{g}_1) \qp
q(g_2,\tilde{g}_2).
\]
As a special case of Proposition \ref{bigth:trivial},
$q(\id_{\group}) = e$. The fact that $q\pair$ depends on $\pair$ 
only up to homotopy is a consequence of the homotopy invariance of
$HF_*\pair$ (Proposition \ref{th:homotopy}) and Theorem \ref{th:gluing}.
\qed

\proof[Proof of Corollary \ref{bigth:invertible}] Because of Lemma
\ref{th:completeness}, this is an immediate consequence of Corollary
\ref{bigth:homomorphism}. \qed

\newcommand{\mymod}{\text{ mod }}
\newcommand{\nequiv}{\not\equiv}
\section{An application to the Maslov index\label{sec:an-application}}

Recall that $\bar{I}: \pi_1(\Ham\mo) \longrightarrow \Z/N\Z$ 
($N$ denotes the minimal Chern number) is defined by
\[
\bar{I}(g) = I\pair \mymod N,
\]
where $\tilde{g}: \tloops \longrightarrow \tloops$ is any lift of
$g \in G$. If $M$ is simply-connected, the definition of
$\bar{I}$ is elementary; in general, it uses Lemma 
\ref{th:contractible-component}, which is based on the Arnol'd conjecture.
Our first result is a simple consequence of the existence of the 
automorphisms $HF_*(g,\tilde{g})$.

\begin{proposition} \label{th:calabi-yau}
If $\mo$ satisfies $c_1|\pi_2(M) = 0$,
$\bar{I}$ is the trivial homomorphism. \end{proposition}

\proof In this case, $\bar{I}$ is $\Z$-valued because
$I\pair \in \Z$ is independent of the choice of $\tilde{g}$.
Since $\bar{I}$ is a homomorphism, it is sufficient to show that
$\bar{I}(g) \leq 0$ for all $g$. The assumption on $c_1$
also implies that the grading of the Novikov ring $\Lambda$ is trivial.
It follows from Theorem \ref{th:what-is-floer-homology}
that $HF_0\mo \iso H_0(M;\Z/2) \otimes \novi \iso \novi$ and
$HF_k\mo = 0$ for $k<0$. By Proposition \ref{th:shift},
$HF_*(g,\tilde{g})$ maps $HF_0\mo$ isomorphically to
$HF_{-2I(g,\tilde{g})}\mo$, which is clearly impossible if
$I(g,\tilde{g}) > 0$. \qed

To obtain more general results, it is necessary to use
the multiplicative structure. Consider
\[
Q^+ = \bigoplus_{i < 2n} H_i(M;\Z/2) \otimes \Lambda \subset \quantum.
\]
Since the degree of any element of $\novi$ is a multiple
of $2N$, $QH_k\mo \subset Q^+$ for any $k$ such that
$k \nequiv 2n \mymod 2N$.

\begin{lemma} \label{th:qplus} If $Q^+ \qp Q^+ \subset Q^+$, every
homogeneous invertible element of $\quantum$ has degree
$2n + 2iN$ for some $i \in \Z$. \end{lemma}

\proof Assume that $x \in QH_k\mo$ is invertible and
$k \nequiv 2n \mymod 2N$; then $x \in Q^+$.
Since $\qp$ has degree $-2n$ and the unit $e$ lies in $QH_{2n}\mo$, 
the inverse $x^{-1}$ has degree $4n-k$. Clearly
$4n - k \nequiv 2n \mymod 2N$, and therefore 
$x^{-1} \in Q^+$. $x^{-1} \qp x = e
\notin Q^+$, contrary to the assumption that $Q^+ \qp
Q^+ \subset Q^+$. \qed

\begin{proposition} If $\mo$ satisfies {\wplus} and
$Q^+ \qp Q^+ \subset Q^+$, $\bar{I}$ is trivial. \end{proposition}

\proof For all $\pair \in \group$, $q\pair \in \quantum$ is
invertible by Corollary \ref{bigth:invertible}.
$q\pair$ has degree $2n - 2I\pair$. Using Lemma \ref{th:qplus}
it follows that $I\pair \equiv 0 \mymod N$, and therefore
$\bar{I}(g) = 0$. \qed

By definition, $Q^+ \qp Q^+ \subset Q^+$ iff
$\gw(x_1,x_2,[pt]) = $ for all $A \in H_2(M;\Z)$
and all $x_1, x_2 \in H_*(M;\Z/2)$ of dimension $<2n$. This is certainly true
if there is an $\o$-tame almost complex structure $J$ on $M$
and a point $y \in M$ such that no non-constant $J$-holomorphic
sphere passes through $y$. A particularly simple case is when there are no
non-constant $J$-holomorphic spheres at all. Then the quantum
intersection product reduces to the ordinary one, that is,
\begin{equation} \label{eq:undeformed}
(x_1 \otimes \gen{\gamma_1}) \qp (x_2 \otimes \gen{\gamma_2}) =
(x_1 \cdot x_2) \otimes \gen{\gamma_1 + \gamma_2}.
\end{equation}
For example, if
\begin{equation} \label{eq:ample-canonical}
(c_1 - \lambda[\o]) |\pi_2(M) = 0
\text{ for some } \lambda<0,
\end{equation}
$c_1(w) < 0$ for any non-constant pseudoholomorphic curve $w$. 
If in addition $N \geq n-2$, there is
a dense set of $J$ such that any $J$-holomorphic sphere satisfies
$c_1(w) \geq 0$ and hence must be constant. 

\begin{corollary} \label{cor:c-negative} If $(M^{2n},\omega)$ 
satisfies \eqref{eq:ample-canonical} and its 
minimal Chern number is at least $n-1$, the homomorphism $\bar{I}$
is trivial. \qed \end{corollary}

Note that we have sharpened the condition on $N$ in order to
fulfil {\wplus}. By a slightly different reasoning,
it seems likely that Corollary \ref{cor:c-negative} 
remains true without any assumption on $N$;
this is one of the motivations for extending the theory
beyond the case where {\wplus} holds. 

As a final example, consider the case where $M$ is four-dimensional.
{\wplus} holds for all symplectic four-manifolds.
For $N = 1$, $\bar{I}$ is vacuous. Assume
that $N \geq 2$ (this implies that $\mo$ is minimal). 
In that case, it is a result of McDuff
that \eqref{eq:undeformed} holds unless $\mo$ is rational or ruled.
For convenience, we reproduce the proof from
\cite{mcduff-salamon96b}:
for generic $J$, there are no non-constant $J$-holomorphic spheres
$w$ with $c_1(w) \leq 0$ because the moduli space of such curves 
has negative dimension. Therefore, if \eqref{eq:undeformed} does
not hold, there is a $w$ with $c_1(w) \geq 2$. After perturbing
$J$, we can assume that $w$ is immersed with transverse
self-intersections \cite[Proposition 1.2]{mcduff91}. Then $\mo$
is rational or ruled by \cite[Theorem 1.4]{mcduff92}.

\begin{corollary} $\bar{I}$ is trivial for all $(M^4,\omega)$
which are not rational or ruled. \qed
\end{corollary}

\newcommand{\Gr}{\mathrm{Gr}}
\section{Examples \label{sec:examples}}

In the first example, $M$ is the Grassmannian $\Gr_k(\C^m)$
with the usual symplectic structure $\omega$ coming from the
Pl{\"u}cker embedding. $\mo$ is a monotone symplectic manifold. 
The $U(m)$-action on $\C^m$ induces a Hamiltonian $PU(m)$-action 
$\rho$ on $M$. Let $g$ be the loop in $\Ham\mo$ given by 
$g_t = \rho(\mathrm{diag}(e^{2\pi i t}, 1, \dots, 1))$. 
Obviously, $[g] \in \pi_1(\Ham\mo)$
lies in the image of $\rho_*: \pi_1(PU(m)) \longrightarrow
\pi_1(\Ham\mo)$, and since $\pi_1(PU(m)) \iso \Z/m$,
$[g]^m = 1$.

Let $H$ be the Hopf bundle over $\CP{1}$.
The fibre bundle $E = E_g$ is the bundle of Grassmannians
associated to the holomorphic vector bundle $H^{-1}
\oplus \C^{m-1}$:
\[
E = \Gr_k(H^{-1} \oplus \C^{m-1}) \stackrel{\pi}{\longrightarrow} \CP{1}.
\]
$H^{-1} \oplus \C^{m-1}$ is a subbundle of the trivial bundle
$\CP{1} \times \C^{m+1}$. This induces a holomorphic map
$f: E \longrightarrow \Gr_k(\C^{m+1})$ which is an embedding on
each fibre $E_z$. Let $\O_z \in \O^2(E_z)$ be the pullback of the
symplectic form on $\Gr_k(\C^{m+1})$ by $f|E_z$. $(E,\O)$ is a
Hamiltonian fibre bundle.

As for any bundle of Grassmannians associated to a vector bundle,
there is a canonical $k$-plane bundle $P_E \longrightarrow E$.
$P_E$ is a subbundle of $\pi^*(H^{-1} \oplus \C^{m-1})$ and
\begin{equation} \label{eq:vertical-tangent-bundle}
TE^v \iso \Hom(P_E,\pi^*(H^{-1} \oplus \C^{m-1})/P_E).
\end{equation}
Therefore $c_1(TE^v) = -k \, c_1(\pi^*H) - m \, c_1(P_E)$.
$P_E$ is isomorphic to $f^*\!P$, where $P$ is the canonical $k$-plane
bundle on $\Gr_k(\C^{m+1})$. It follows that
\begin{equation} \label{eq:c-one}
\deg(s^*TE^v) = -k -m \deg( s^*\!f^*\!P)
\end{equation}
for any section $s$ of $E$. If $s$ is holomorphic, $f(s)$ is a 
holomorphic curve in $\Gr_k(\C^{m+1})$, and since 
$c_1(P) \in H^2(\Gr_k(\C^{m+1}))$ is a negative multiple of the
symplectic class, either $\deg(s^*\!f^*\!P) \leq -1$ or $f(s)$ is constant.
In the first case, $\deg(s^*TE^v) > 0$ by \eqref{eq:c-one}. In the
second case, $s$ must be one of the `constant' sections
\[
s_W(z) = 0 \oplus W \in \Gr_k(H_z \oplus \C^{m-1})
\]
for $W \in \Gr_k(\C^{m-1})$.
We now check that $(E,\O)$ satisfies the conditions
of Proposition \ref{th:computation}.

\begin{enumerate} \item We have shown that any holomorphic section
with $\deg(s^*TE^v) \leq 0$ is a `constant' one. The space $\sec$
of such sections is certainly connected. Their $\Gamma$-equivalence
class $S_0$ satisfies $c_1(TE^v)(S_0) = -k$.
\item $s_W^*TE^v = \Hom(W, H^{-1} \oplus \C^{m-1}/W) \iso
\Hom(\C^k, H^{-1} \oplus \C^{m-k-1})$ is a sum of line bundles
of degree $0$ or $-1$, hence $H^{0,1}(\CP{1},s_W^*TE^v) = 0$.
\item and (iv) follow from the fact that $\mo$ is monotone with minimal
Chern number $N = m > k$. \end{enumerate}

$(E_{z_0},\O_{z_0})$ is identified with $\mo$ by choosing an
element of $H_{z_0}^{-1}$ which has length $1$ for the standard Hermitian
metric on $H^{-1}$. The evaluation map $\ev_{z_0}: \sec \longrightarrow M$
is an embedding whose image is $\Gr_k(\C^{m-1}) \subset \Gr_k(\C^m)$.
By Proposition \ref{th:computation},
$Q(E,\O,S_0) = [\Gr_k(\C^{m-1})] \otimes \gen{0}$.
Using the diagram \eqref{diag:basic}, we can derive from $[g]^m = 1$ that 
\begin{equation} \label{eq:relation}
\left([\Gr_k(\C^{m-1})] \otimes \gen{0}\right)^m = [M] \otimes \gen{\gamma}
\end{equation}
for some $\gamma \in \Gamma$. Because of the
grading, $\gamma$ must be $k$ times the standard generator
of $\pi_2(M) \iso H_2(M;\Z)$. \eqref{eq:relation} can be verified
by a direct computation, since $(QH_*\mo,\qp)$ is known
\cite{siebert-tian94} \cite{witten93}. \vspace{0.5em}

Our second example concerns the rational ruled surface
$M = \ruled_2$. Recall that $\ruled_r$ $(r>0)$
is the total space of the holomorphic fibre bundle
$p_r: \proj(\C \oplus H^r) \longrightarrow \CP{1}$. We will use
two standard facts about $\ruled_r$.

\begin{lemma} \label{th:beauville} For any holomorphic
section $s$ of $p_r$, 
$
\deg(s^*T\ruled_r^v) \geq -r.
$
There is a unique section such that equality holds, and all
others have $\deg(s^*T\ruled_r^v) \geq -r+2$. \end{lemma}

\begin{lemma} \label{th:beauville-b} For every
holomorphic map $w: \CP{1} \longrightarrow \ruled_2$,
$c_1(w) \geq 0$. All non-constant $w$ with $c_1(w) < 2$ are
of the form $w = s_- \circ u$, where $u: \CP{1} \longrightarrow \CP{1}$
is holomorphic and $s_-$ is the section
with $\deg(s_-^*T\ruled_2^v) = -2$. \end{lemma}

Lemma \ref{th:beauville} is a special case of a theorem on
irreducible curves in $\ruled_r$ \cite[Proposition IV.1]{beauville}. Lemma
\ref{th:beauville-b} can be derived from Lemma \ref{th:beauville}
by considering a map $w: \CP{1} \longrightarrow \ruled_2$ 
such that $p_2w$ is not constant as a
section of the pullback $(p_2w)^*\ruled_2$.

Since $\proj(\C \oplus H^2) \iso \proj(H^{-2} \oplus \C)$, an embedding
of $H^{-2} \oplus \C$ into the trivial bundle $\CP{1} \times \C^3$ defines a
holomorphic map $f: M \longrightarrow \CP{2}$. Let $\tau_k$
be the standard integral K{\"a}hler form on $\CP{k}$. For all
$\lambda>1$, $\o_\lambda = (\lambda-1) p_2^*\tau_1 + f^*\tau_2$ is a
K{\"a}hler form on $M$.

The action of $S^1$ on $H^2$ by multiplication induces a Hamiltonian
circle action $g$ on $(M,\o_\lambda)$. $(E,\O) = (E_g,\O_g)$
can be constructed as follows: $E = \proj(V)$ is the bundle of projective
spaces associated to the holomorphic vector bundle
\[
V = \C \oplus \pr_1^*H^2 \otimes \pr_2^*H^{-1}
\]
over $\CP{1} \times \CP{1}$ ($\pr_1,\pr_2$ are
the projections from $\CP{1} \times \CP{1}$ to $\CP{1}$).
$\pi: E \longrightarrow \CP{1}$  is obtained by composing
$\pi_V: \proj(V) \longrightarrow \CP{1} \times \CP{1}$ with $\pr_2$.
The fibres $E_z = \pi^{-1}(z)$ are the ruled surfaces
$\proj(\C \oplus H^2 \otimes H^{-1}_z)$. Any $\xi \in H^{-1}_z$ 
with unit length for the standard Hermitian metric determines a biholomorphic
map $E_z \longrightarrow M$. Since two such maps differ by
an isometry of $(M,\o_\lambda)$, the symplectic structure $\O_z$
on $E_z$ obtained in this way is independent of the choice of $\xi$.
$(E,\O)$ is a Hamiltonian fibre bundle.

A holomorphic section $s$ of $\pi$ can be decomposed into two pieces:
a section $s_1 = \pi_V \circ s: \CP{1} \longrightarrow \CP{1} \times
\CP{1}$ of $\pr_2$ and a section $s_2$ of
\[
F = \proj(s_1^*V) \longrightarrow \CP{1}.
\]
In a sense, this second piece is given by $s$ itself, using the fact that
$s(z) \in \proj(V_{s_1(z)})$ for all $z$. The decomposition leads to
an exact sequence
\begin{equation} \label{eq:vector-bundle-sequence}
0 \longrightarrow s_2^*(TF^v) \longrightarrow
s^*TE^v \stackrel{D\pi_V}{\longrightarrow} s_1^*(\ker(D\pr_2))
\longrightarrow 0
\end{equation}
of holomorphic vector bundles over $\CP{1}$. Since $s_1$ is a 
section of $\pr_2$, it is given by $s_1(z) = (u(z),z)$ for some
$u: \CP{1} \longrightarrow \CP{1}$. Clearly, 
$s_1^*(\ker(D\pr_2)) = u^*T\CP{1}$. Let $d$ be the degree of
$u$. The fibre bundle $F = \proj(\C \oplus u^*H^2 \otimes
H^{-1})$ is isomorphic to $\ruled_1$ for $d = 0$ and to
$\ruled_{2d-1}$ for $d > 0$. If $d>0$, 
$\deg(s_2^*TF^v) \geq 1-2d$ by Lemma \ref{th:beauville}. Since
$\deg(u^*T\CP{1}) = 2d$, it follows from 
\eqref{eq:vector-bundle-sequence} that $\deg(s^*TE^v) > 0$.
If $d = 0$, $s_1(z) = (c,z)$ for some $c \in \CP{1}$. The same 
argument as before shows that $\deg(s^*TE^v) > 0$ unless $s_2$ 
is the unique section of $F$ such that $\deg(s_2^*TF^v) = -1$. 
It is easy to write down this section explicitly. We conclude
that any holomorphic section $s$ of $\pi$ with 
$\deg(s^*TE^v) \leq 0$ belongs to the family
$\sec = \{s_c\}_{c \in \CP{1}}$, where
\begin{equation} \label{eq:explicit-section}
s_c(z) = [1:0] \in \proj(\C \oplus H_c^2 \otimes H_z^{-1}) \subset E_z.
\end{equation}
These sections have $\deg(s_c^*TE^v) = -1$.
For any $z \in \CP{1}$, the evaluation map $\sec \longrightarrow E_z$
is an embedding. Its image is the curve of self-intersection $2$
which corresponds to $C^+ = \proj(\C \,\oplus\, 0) \subset M$ under
an isomorphism $E_z \iso M$ as above. 
We will now verify that the conditions of Proposition
\ref{th:computation} are satisfied. 

\begin{enumerate} \item We have already shown that
$\sec$ is connected and that $c_1(TE^v)(S_0) = -1$ for
its $\Gamma$-equivalence class $S_0$.
\item For $s \in \sec$, \eqref{eq:vector-bundle-sequence} reduces to
\[
0 \longrightarrow H^{-1} \longrightarrow s^*TE^v \longrightarrow \C
\longrightarrow 0.
\]
By the exact sequence of cohomology groups, 
$H^{0,1}(\CP{1},s^*TE^v) = 0$.
\item is part of Lemma \ref{th:beauville-b}.
\item Consider the curve $C^- = \proj(0 \oplus H^2) \subset M$.
$C^-$ is the image of the unique section $s_-$ of $p_2$ with
$\deg(s_-^*T\ruled_2^v) = -2$. Let $w: \CP{1} \longrightarrow E_z$
be a non-constant holomorphic map with $c_1(TE)(w) < 2$. 
By Lemma \ref{th:beauville-b}, the image of $w$ is mapped
to $C^-$ under an isomorphism $i: E_z \iso M$ chosen as above. We have
seen that $i(s(z)) \in C^+$ for all $s \in \sec$.
Since $C^+ \cap C^- = \emptyset$, it follows that $s(z) \notin \im(w)$.
\end{enumerate}

We obtain $Q(E,\O,S_0) = [C^+] \otimes \gen{0}$.  

\begin{lemma} \label{th:quantum-relation} 
Let $x^+, x^- \in \pi_2(M) \iso H_2(M;\Z)$ be the classes of
$C^+, C^-$, and $\bar{x}^+,\bar{x}^- \in H_2(M;\Z/2)$ their
reductions mod $2$. Then
\[
(\bar{x}^+ \otimes \gen{0})^2 =
[M] \otimes \left(\gen{\tfrac{1}{2}(x^+ - x^-)} -
\gen{\tfrac{1}{2}(x^+ + x^-)}\right)
\]
in $(QH_*(M,\o_\lambda),\qp)$. \end{lemma}

\proof In \cite{mcduff87} McDuff showed that $(M,\o_\lambda)$ is
symplectically isomorphic to $\CP{1} \times \CP{1}$ with the
product structure $\lambda(\tau_1 \times 1) + 1 \times \tau_1$.
Such an isomorphism maps $x^{\pm}$ to
$a \pm b$, where $a = [\CP{1} \times \mathit{pt}]$ and
$b = [\mathit{pt} \times \CP{1}]$. Let $\bar{a}, \bar{b}$ the mod $2$
reductions of these classes. The
quantum intersection product on $\CP{1} \times \CP{1}$ is known
(see e.g. \cite[Proposition 8.2 and Example 8.5]{ruan-tian94}); it satisfies
\[
(\bar{a} \otimes \gen{0})^2 = [\CP{1} \times \CP{1}] \otimes \gen{b}, \;
(\bar{b} \otimes \gen{0})^2 = [\CP{1} \times \CP{1}] \otimes \gen{a}.
\]
Because of the $\Z/2$-coefficients, this implies the relation stated
above. \qed

This can be used to give a proof of the following result of McDuff.

\begin{cor} For all $\lambda > 1$, 
$[g] \in \pi_1(\Ham(M,\omega_\lambda))$ has infinite order. \end{cor}

\proof Let $\tilde{g}: \tloops \rightarrow \tloops$ be the lift
of $g$ corresponding to the equivalence class $S_0$. 
By Theorem \ref{bigth:homomorphism} and Lemma
\ref{th:quantum-relation},
\[
q(g^2,\tilde{g}^2) = Q(E,\O,S_0)^2 =
[M] \otimes \gen{\tfrac{1}{2}(x^+ - x^-)}(\gen{0} - \gen{x^-})
\]
and
\begin{align*}
q(g^{2m},\tilde{g}^{2m}) &= [M] \otimes \gen{\tfrac{m}{2}(x^+ - x^-)}
(\gen{0} - \gen{x^-})^m,\\
q(g^{2m+1},\tilde{g}^{2m+1}) &= \bar{x}^+ \otimes
\gen{\tfrac{m}{2}(x^+ - x^-)} (\gen{0} - \gen{x^-})^m
\end{align*}
for all $m \geq 0$. Since $x^-$ is represented by an algebraic curve,
$\o_\lambda(x^-) > 0$ for all $\lambda > 1$ (in fact, 
$\o_\lambda(x^-) = \lambda - 1$). Therefore the class of
$x^-$ in $\Gamma$ has infinite order, and
\[
(\gen{0} - \gen{x^-})^m 
\notin \{\gen{\gamma} \; | \; \gamma \in \Gamma\} \subset \Lambda
\]
for all $m \geq 1$. It follows that $q(g^k,\tilde{g}^k) \notin
\tau(\Gamma)$ for all $k>0$. Because of the diagram \eqref{diag:basic},
this implies that $[g^k] \in \pi_1(\Ham\mo)$ is nontrivial. \qed

\begin{rem} Using again Theorem \ref{bigth:homomorphism} and
Lemma \ref{th:quantum-relation}, one computes that
\begin{multline*}
q(g^{-1},\tilde{g}^{-1}) = (\bar{x}^+ \otimes \gen{0})^{-1} =\\
= \bar{x}^+ \otimes \gen{\tfrac{1}{2}(x^- - x^+)}
\left(\gen{0} + \gen{x^-} + \gen{2x^-} + \cdots \right).
\end{multline*}
In this case, a direct computation of the invariant seems to be more 
difficult because infinitely many moduli spaces contribute to it.
\end{rem}

\bibliographystyle{amsplain}

\begin{thebibliography}{10}

\bibitem{abreu97}
M.~Abreu, \emph{Topology of symplectomorphism groups of {$S^2 \times S^2$}},
  Preprint.

\bibitem{banyaga78}
A.~Banyaga, \emph{Sur la structure du groupe des diff{\'e}omorphismes qui
  preservent une forme symplectique}, Commun. Math. Helv. \textbf{53} (1978),
  174--227.

\bibitem{beauville}
A.~Beauville, \emph{Complex algebraic surfaces}, Cambridge Univ. Press, 1983.

\bibitem{floer88}
A.~Floer, \emph{Symplectic fixed points and holomorphic spheres}, Commun. Math.
  Phys. \textbf{120} (1989), 575--611.

\bibitem{floer-hofer93}
A.~Floer and H.~Hofer, \emph{Coherent orientations for periodic orbit problems
  in symplectic geometry}, Math. Z. \textbf{212} (1993), 13--38.

\bibitem{gromov85}
M.~Gromov, \emph{Pseudoholomorphic curves in symplectic manifolds}, Invent.
  Math. \textbf{82} (1985), 307--347.

\bibitem{hofer-salamon95}
H.~Hofer and D.~Salamon, \emph{Floer homology and {N}ovikov rings}, The {F}loer
  memorial volume (H.~Hofer, C.~Taubes, A.~Weinstein, and E.~Zehnder, eds.),
  Progress in Mathematics, vol. 133, Birkh{\"a}user, 1995, pp.~483--524.

\bibitem{mcduff87}
D.~McDuff, \emph{Examples of symplectic structures}, Invent. Math. \textbf{89}
  (1987), 13--36.

\bibitem{mcduff91}
\bysame, \emph{The local behaviour of holomorphic curves in almost complex
  $4$-manifolds}, J. Differential Geom. \textbf{34} (1991), 143--164.

\bibitem{mcduff92}
\bysame, \emph{Immersed spheres in symplectic $4$-manifolds}, Ann. Inst.
  Fourier \textbf{42} (1992), 369--392.

\bibitem{mcduff-salamon}
D.~McDuff and D.~Salamon, \emph{{$J$}-holomorphic curves and quantum
  cohomology}, University Lecture Notes Series, vol.~6, Amer. Math. Soc., 1994.

\bibitem{mcduff-salamon96}
\bysame, \emph{Introduction to {S}ymplectic {T}opology}, Oxford University
  Press, 1995.

\bibitem{mcduff-salamon96b}
\bysame, \emph{A survey of symplectic $4$-manifolds with $b^+ = 1$}, Turkish J.
  Math. \textbf{20} (1996), 47--60.

\bibitem{ono94}
K.~Ono, \emph{On the {A}rnold conjecture for weakly monotone symplectic
  manifolds}, Invent. Math. \textbf{119} (1995), 519--537.

\bibitem{piunikhin-salamon-schwarz94}
S.~Piunikhin, D.~Salamon, and M.~Schwarz, \emph{Symplectic {F}loer-{D}onaldson
  theory and quantum cohomology}, Contact and symplectic geometry (C.~B.
  Thomas, ed.), Cambridge Univ. Press, 1996, pp.~171--200.

\bibitem{ruan-tian94}
Y.~Ruan and G.~Tian, \emph{A mathematical theory of {Q}uantum {C}ohomology}, J.
  Differential Geom. \textbf{42} (1995), 259--367.

\bibitem{salamon-zehnder92}
D.~Salamon and E.~Zehnder, \emph{Morse theory for periodic solutions of
  {H}amiltonian systems and the {M}aslov index}, Comm. Pure Appl. Math.
  \textbf{45} (1992), 1303--1360.

\bibitem{schwarz97}
M.~Schwarz, \emph{An explicit isomorphism between {F}loer homology and quantum
  cohomology}, In preparation.

\bibitem{schwarz96}
\bysame, \emph{A quantum cup-length estimate for symplectic fixed points},
  Preprint.

\bibitem{schwarz}
\bysame, \emph{Morse homology}, Progress in Mathematics, vol. 111,
  Birkh{\"a}user, 1993.

\bibitem{schwarz95}
\bysame, \emph{Cohomology operations from {$S^1$}-cobordisms in {F}loer
  homology}, Ph.D. thesis, {ETH} {Z}{\"u}rich, 1995.

\bibitem{siebert-tian94}
B.~Siebert and G.~Tian, \emph{On quantum cohomology rings of {F}ano manifolds
  and a formula of {V}afa and {I}ntriligator}, Preprint.

\bibitem{thurston76}
W.~Thurston, \emph{Some simple examples of symplectic manifolds}, Proceedings
  of the Amer. Math. Soc. \textbf{55} (1976), 467--468.

\bibitem{weinstein89}
A.~Weinstein, \emph{Cohomology of symplectomorphism groups and critical values
  of {H}amiltonians}, Math. Z. \textbf{201} (1989), 75--82.

\bibitem{witten88}
E.~Witten, \emph{Topological {S}igma models}, Commun. Math. Phys. \textbf{118}
  (1988), 411--449.

\bibitem{witten93}
\bysame, \emph{The {V}erlinde algebra and the cohomology of the
  {G}rassmannian}, Geometry, topology and physics (S.-T. Yau, ed.),
  International Press, 1995, pp.~357--422.

\end{thebibliography}
\providecommand{\bysame}{\leavevmode\hbox to3em{\hrulefill}\thinspace}

\end{document}